\documentclass[a4paper,12pt]{elsarticle}
\journal{M2AN}
\usepackage{amsmath}
\usepackage[english]{babel}
\usepackage{graphicx}
\usepackage[usenames,dvipsnames]{color} 
\usepackage[numbers]{natbib}
\usepackage[footnotesize,bf,sf,center]{caption}
\usepackage{float}
\usepackage[dvips]{epsfig}
\usepackage[bottom=2.5cm,top=2.5cm,left=3cm,right=2cm]{geometry}
\usepackage{calc}
\usepackage{soul}
\usepackage{amsmath,amsfonts,amssymb}
\usepackage{subfigure}
\usepackage{setspace}
\usepackage{amssymb}
\usepackage{mathrsfs}
\usepackage{algorithm}
\usepackage{algpseudocode}
\usepackage{dsfont}
\usepackage{epstopdf}
\usepackage{subfig}%
\usepackage[utf8]{inputenc}
\usepackage{nicefrac,xfrac}
\newcommand{\ve}[1]{\mbox{\boldmath $#1$}}
\newcommand{\vesmall}[1]{\mbox{\boldmath \scriptsize $#1$}}

\newtheorem{theorem}{Theorem}[section]

\newtheorem{lemma}[theorem]{Lemma}

\newtheorem{remark}{Remark}[section]

\newdefinition{rmk}{Remark}
\newcommand{\proof} [1]{ \noindent {\bf Proof.} #1 \hfill\rule{0.5em}{1.2ex} \par\medskip}
\usepackage{graphicx}
\begin{document}
\begin{frontmatter}

  \renewcommand\arraystretch{1.0}

    \title{\textbf{A fully segregated and unconditionally stable IMEX scheme for dispersed multiphase flows}}
    \author{
   {\bf Douglas R.~Q.~Pacheco} $^{1,2,3}$ and \ 
   {\bf Richard Schussnig} $^{4}$, 
     \\
   {\small ${}^{1}$ Chair for Computational Analysis of Technical Systems, RWTH Aachen University, Germany}\\
   {\small ${}^{2}$ Chair of Methods for Model-based Development in Computational Engineering, RWTH}\\
   {\small ${}^{3}$ Center for Simulation and Data Science (JARA-CSD), RWTH Aachen University, Germany}\\
   {\small ${}^{4}$ Faculty of Mathematics, Ruhr University Bochum, Germany}}
\begin{keyword}
Multi-phase flow \sep Two-fluid system \sep Fractional-step methods \sep Projection schemes \sep Euler-Euler system  \sep Eulerian-Eulerian model
\end{keyword}
\begin{abstract}
Euler--Euler or volume-averaged Navier--Stokes equations are used in various applications to model systems with two or more interpenetrating phases. Each fluid obeys its own momentum and mass equations, and the phases are typically coupled via drag forces and a shared pressure. Monolithic solvers can therefore be very expensive and difficult to implement, so there is great computational appeal for decoupled methods. However, splitting the subproblems requires treating the coupling terms (pressure and drag) explicitly, which must be done carefully to avoid time-step restrictions. In this context, we derive a new first-order pressure-correction method based on the incompressibility of the mean velocity field, combined with an explicit treatment of the drag forces. Furthermore, both the convective and viscous terms are treated semi-implicitly. This gives us an implicit-explicit (IMEX) method that is very robust not only due to its unconditional energy stability, but also because it does not require any type of fixed-point iterations. Each time step has only linear, scalar transport equations and a single pressure Poisson problem as building blocks. We rigorously prove temporal stability without any CFL-like conditions, and the theory is confirmed through two-phase numerical examples.
\end{abstract}
\end{frontmatter}
	
\section{Introduction}
    
    In Eulerian approaches, a dispersed multiphase system is typically modelled as a set of two or more fluids that completely share the same domain $\Omega$, with volume fractions describing how much of each fluid occupies each point in $\Omega$ at a given time. These models are different from \textsl{separated} two-phase flows, where each fluid occupies its own \textsl{subset} of $\Omega$, separated by either sharp \cite{Gross2013} or diffuse interfaces, as in phase-field approaches \cite{Shen2010,Deteix2022,Salgado2013,Liu2015,AcostaSoba2025}---refer to the recent article \cite{GarciaVillalba2025} for an overview of different multiphase models. In the dispersed or volume-averaged setting, which is the topic of this article, how to model the interactions between the phases depends, e.g., on the flow regime and the materials considered \cite{Caia2004,Varaksin2013}. These models are relevant in various industrial and scientific applications, such as gas-liquid systems \cite{Triplett1999,Percival2014}, debris flows \cite{Rauter2021,Rauter2022} and volcanic ash dynamics \cite{Jacobs2013}. Regardless of how the phase interactions are described, most volume-averaged models have in common that (1) each fluid obeys its own Navier--Stokes-like system in $\Omega$, (2) the volume fractions form (together) a partition of unity, and (3) there is a pressure field at least partially shared by all the phases. 
    
    Regarding pressure modelling, the most popular approach in incompressible flows---notably in environmental applications \cite{Jacobs2013,Kazeminezhad2012}---is to consider a single pressure common to all phases. In some scenarios such as dense granular flows, a so-called solid pressure is typically added for the granular phase, but this pressure is modelled as a function of other flow quantities instead of being treated as an additional unknown \cite{Rauter2021,Rauter2022,Zhang2021}. Multi-pressure models also exist, where each phase has its own pressure, all of  which being additional unknowns. These models are particularly important in compressible gas dynamics and reacting flows to account for pressure non-equilibrium effects \cite{Baer1986,Saurel1999}. In this work, we consider incompressible regimes, focusing on the prototypical variant where the phases are coupled by a common pressure, in addition to drag forces \cite{Gravenkamp2024,Pacheco2025}.
    
	The numerics of dispersed multiphase flows pose considerable challenges. Since each phase follows a Navier--Stokes-like system, the well-known issues and instabilities observed in single-phase flow simulation also arise. Furthermore, the number of unknowns can be substantial: an $M$-phase flow will involve $M$ phase fractions, $M$ velocities and a pressure. Probably due to the complexity of such systems, numerical methods designed specifically for them are scarce in the literature. In the past thirty years, some works have emerged addressing different types of numerical issues, such as absent volume-fraction diffusion \cite{Hiltunen1997,Vreman2011,Dang2014} or poor mass conservation \cite{Gravenkamp2024}. Attempts have also been made towards decoupling the phases via fractional stepping \cite{Caia2004,Jacobs2013,Behrangi2019}, but until very recently nothing beyond purely heuristic approaches. In fact, numerical stability is such a delicate matter in multiphase flows, that even schemes with implicit pressures often suffer from \textsl{parabolic} CFL restrictions \cite{Rauter2021,Ferrari2024}. To the best of our knowledge, the first provenly stable fractional-step scheme for dispersed multiphase flows was just recently proposed \cite{Pacheco2025}. Although that method decouples the phases while still being energy stable, it introduces an auxiliary pressure for each phase, which increases the number of unknowns and becomes increasingly less efficient as the number of phases grows.  
    
    This work improves our previous approach \cite{Pacheco2025} in three important ways: (1) the new scheme has only one pressure Poisson equation regardless of the number of phases, (2) we devise a simpler, fully explicit treatment of the drag term, and (3) our new implicit-explicit (IMEX) treatment of the viscous term allows splitting even each velocity update into a series of scalar subproblems. IMEX methods are used in many applications to circumvent the cost and complexity of fully implicit schemes, while still avoiding the severe stability restrictions of fully explicit methods. In incompressible flows, convective and stabilisation terms may be treated (semi-)explicitly to avoid nonlinearities and reduce computational costs \cite{Burman2023,Burman2024}, while viscous terms must usually be kept implicit \cite{Barrenechea2024}. In compressible gas dynamics, on the other hand, severe time-step restrictions can be avoided through IMEX methods that treat acoustic terms implicitly \cite{Boscarino2019}.

    The scheme presented here can be seen as an (IMEX) extension of the popular incremental pressure-correction method \cite{Guermond2006} to multiphase flows. As a matter of fact, our inspiration comes from different families of fractional-step methods for incompressible flows \cite{Deteix2022,Guermond2009}. The resulting scheme allows using single-phase Navier--Stokes solvers---or even \textsl{scalar} advection-diffusion-reaction solvers---as building blocks for the multiphase system. Moreover, each subproblem becomes linearised, which further simplifies implementation and improves both efficiency and robustness (the issue of the often delicate convergence of Newton or Picard schemes is eliminated). These advantages are attained without sacrificing numerical stability, which we rigorously prove to be unconditional with respect to the time-step size. While the focus here is on the temporal discretisation, our scheme should be relatively simple to implement within various spatial frameworks, especially finite element-based ones. 
	
    We organise the rest of this article as follows. Section \ref{sec_pre} introduces the model problem and a consistent reformulation aimed at improved numerical stability in the discrete case. The conservation properties of the continuous system are derived and discussed in Section \ref{sec_conservation}. The fractional-step method is derived in Section \ref{sec_schemes}, where consistency aspects are also discussed. The unconditional stability is then proved in Section \ref{sec_analysis} for pressure and velocities. Section \ref{sec_convection} addresses stability issues related to the hyperbolic volume-fraction problem and proposes some remedies. In Section \ref{sec_examples} we provide numerical examples assessing the accuracy and the stability of our new method, before drawing concluding remarks in Section \ref{sec_Conclusion}.
	
	\section{Preliminaries}\label{sec_pre}
	\subsection{Model problem}
	The equations tackled here can be derived by volume-averaging the Navier--Stokes equations for each phase and making certain simplifying and modelling assumptions. While the derivation is out of our scope, details can be found in the literature \cite{Caia2004,Gravenkamp2024,Rusche2003}. We consider standard notation for Hilbert spaces. The $L^2(\Omega)$ product of two functions is represented by $\langle \cdot ,\cdot \rangle$, with $\| \cdot \|$ and $\| \cdot \|_{\infty}$ denoting the $L^2(\Omega)$ and $L^{\infty}(\Omega)$ norms. The data, including initial and boundary conditions, are assumed as sufficiently smooth. 
    
    For a time interval $(0,T]$ and a domain $\Omega\subset\mathbb{R}^{d}$, $d=2$ or $3$, the balance equations for the flow of the $k$th phase in $Q:=\Omega\times(0,T]$ are
	\begin{flalign}
		&\partial_t\alpha_k + \nabla\cdot(\alpha_k\ve{u}_k)  = 0\, , \label{mass}\\
		&\rho_k[\partial_t(\alpha_k\ve{u}_k) + \nabla\cdot(\alpha_k\ve{u}_k\otimes\ve{u}_k - 2\alpha_k\nu_k\nabla^{\mathrm{s}}\ve{u}_k)] + \alpha_k\nabla p + \sum_{l=1}^{M}\gamma_{kl}(\ve{u}_k-\ve{u}_l) = \rho_k\alpha_k\ve{g}_k\, , \label{momentum}
	\end{flalign}
for $k=1,\dots,M$, with all the $M$ phases occupying $Q$ such that
\begin{align}
  \sum_{k=1}^M\alpha_k \equiv 1\, . \label{sum1}
\end{align}
The unknowns are a single pressure $p$, $M$ phase velocities $\ve{u}_k$, and $M$ volume fractions $\alpha_k$. The symmetric velocity gradient is defined in the usual way:
\begin{align*}
    \nabla^{\mathrm{s}}\ve{u} = \frac{1}{2}\big(\nabla\ve{u} + \nabla^{\top}\ve{u}\big)\, ,
\end{align*}
where $\nabla^{\top}(\cdot) = [\nabla (\cdot)]^{\top}$ denotes the transpose of the gradient. The model parameters for each phase are a kinematic viscosity $\nu_k>0$, a (constant) density $\rho_k$, a volumetric force $\ve{g}_k$, and the drag coefficients $\gamma_{kl}\geq 0$, which may depend on a combination of $|\ve{u}_k-\ve{u}_l|$, $\alpha_k$ and $\alpha_l$.  Because $\gamma_{kk} = 0$ and $\gamma_{kl} = \gamma_{lk}$, only $\nicefrac12 M(M-1)$ (not $M^2$) drag coefficients are needed.  have $\ve{u}_k$ defined everywhere in $\Omega$, each $\alpha_k$ must be strictly positive; the presence of regions with one or more null volume fractions requires considering moving fluid domains (i.e., free-surface problems), which can be done more naturally within arbitrary Lagrangian--Eulerian (ALE) frameworks \cite{Dang2014}.
	
\subsection{Equivalent reformulation}
An essential ingredient for the stability of the numerical scheme presented later is a consistent reformulation of the governing equations. We start by multiplying Eq.~\eqref{mass} by $-\frac12\rho_k\ve{u}_k$ and adding the result to Eq.~\eqref{momentum}. Then, denoting $\lambda_k := \sqrt{\rho_k\alpha_k}$, the inertial terms $\rho_k\left[\partial_t(\alpha_k\ve{u}_k) + \nabla\cdot(\alpha_k\ve{u}_k\otimes\ve{u}_k)\right]$ in the momentum equations can be rewritten as
\begin{align*}	
&\rho_k\left[\partial_t(\alpha_k\ve{u}_k) + \nabla\cdot(\alpha_k\ve{u}_k\otimes\ve{u}_k)\right]\\
&=\rho_k\bigg\lbrace\partial_t(\alpha_k\ve{u}_k) + \nabla\cdot(\alpha_k\ve{u}_k\otimes\ve{u}_k)  - \frac{1}{2}[\overbrace{\partial_t\alpha_k + \nabla\cdot(\alpha_k\ve{u}_k)}^{=\, 0 \, \text{(consistent with \eqref{mass})}}]\ve{u}_k\bigg\rbrace \\
        &= \rho_k\left\lbrace\partial_t(\alpha_k\ve{u}_k) + \alpha_k\ve{u}_k\cdot\nabla\ve{u}_k + [\nabla\cdot(\alpha_k\ve{u}_k)]\ve{u}_k  - \frac{1}{2}[\partial_t\alpha_k + \nabla\cdot(\alpha_k\ve{u}_k)]\ve{u}_k\right\rbrace\\
        &= \partial_t(\lambda_k\lambda_k\ve{u}_k)  + \lambda_k^2\ve{u}_k\cdot\nabla\ve{u}_k + \frac{1}{2}[\nabla\cdot(\lambda_k^2\ve{u}_k)]\ve{u}_k - \frac{1}{2}[\partial_t(\lambda_k^2)]\ve{u}_k   \\
        &= \lambda_k\partial_t(\lambda_k\ve{u}_k) + (\partial_t\lambda_k)(\lambda_k\ve{u}_k) + \lambda_k^2\ve{u}_k\cdot\nabla\ve{u}_k + \frac{1}{2}[\nabla\cdot(\lambda_k^2\ve{u}_k)]\ve{u}_k - (\lambda_k\partial_t\lambda_k)\ve{u}_k\\
        &= \rho_k\left\lbrace\sqrt{\alpha_k}\partial_t(\sqrt{\alpha_k}\ve{u}_k) + \alpha_k\ve{u}_k\cdot\nabla\ve{u}_k + \frac{1}{2}[\nabla\cdot(\alpha_k\ve{u}_k)]\ve{u}_k  \right\rbrace,
\end{align*}
with the notation $\ve{v}\cdot\nabla\ve{u} := (\nabla\ve{u})\ve{v}$. This ``conservative'' form is important to guarantee temporal stability for the velocities, as will be shown in Section \ref{sec_analysis}.
   
For our next reformulation, we add up Eqs.~\eqref{mass} for all $M$ phases:  
	\begin{align*}
		\sum_{k=1}^{M}\nabla\cdot(\alpha_k\ve{u}_k) &= -\partial_t \sum_{k=1}^M\alpha_k\\
		&= -\partial_t(1)\\
		&=0\, ,
	\end{align*}
	thanks to \eqref{sum1}. This divergence-free condition on the mean velocity $\bar{\ve{u}} := \sum_{k=1}^{M}\alpha_k\ve{u}_k$ will be used instead of the algebraic constraint \eqref{sum1}.
	
Finally, denoting $\mu_k:=\rho_k\nu_k$, our \textsl{consistently} modified system reads
	\begin{align}
		\partial_t\alpha_k + \ve{u}_k\cdot\nabla\alpha_k + (\nabla\cdot\ve{u}_k)\alpha_k &= 0\, ,\quad  k=1,\dots,M,  \label{mass_}\\
		\rho_k\bigg[\sqrt{\alpha_k}
    \partial_t(\sqrt{\alpha_k}\ve{u}_k) + \alpha_k\ve{u}_k\cdot\nabla\ve{u}_k  + \frac{1}{2} \nabla\cdot(\alpha_k\ve{u}_k)\ve{u}_k\bigg]&- \nabla\cdot(2\alpha_k\mu_k\nabla^{\mathrm{s}}\ve{u}_k)\nonumber \\
    +\, \alpha_k\nabla p   + \sum_{l=1}^{M}\gamma_{kl}(\ve{u}_k-\ve{u}_l)&= \rho_k\alpha_k\ve{g}_k\, , \quad  k=1,\dots,M ,  \label{momentum2}\\
\nabla\cdot\sum_{k=1}^M\alpha_k\ve{u}_k &= 0\, ,\label{divFree}
	\end{align}
equipped with appropriate initial and boundary conditions. Notably, we need the initial volume fractions $\alpha_k(0)$ to be such that  $\sum_{k=1}^M\alpha_k(0) \equiv 1$, since combining \eqref{mass_} and \eqref{divFree} implies
\begin{align*}
    \partial_t\left(\sum_{k=1}^M\alpha_k\right) = 0\, ,
\end{align*}
which will then give us
\begin{align*}
    \sum_{k=1}^M\alpha_k(t) = \sum_{k=1}^M\alpha_k(0) \equiv 1\, ,
\end{align*}
thereby guaranteeing the equivalence.
\begin{remark}
   Alternatively, we could replace one (any) of the $M$ equations \eqref{mass_} by the unity constraint \eqref{sum1} and still maintain the equivalence. What is most important for the decoupling scheme introduced later is that the divergence-free constraint \eqref{divFree} be one of the equations in the (continuous) system.
\end{remark}

\section{Conservation properties}\label{sec_conservation}
\subsection{Energy}
The reformulation above makes the momentum equations convenient for establishing an energy balance. To do that, we will consider a scenario without external energy sources ($\ve{g}_k=\ve{0}$ and $\ve{u}_k|_{\partial\Omega} = \ve{0}$). Dotting \eqref{momentum2} with $\ve{u}_k$ and integrating over $\Omega$ gives us
\begin{align*}
    A_k + B_k + C_k + D_k + E_k = 0\, ,
\end{align*}
where
\begin{align*}
  &A_k = \rho_k\int_{\Omega}(\sqrt{\alpha_k}\ve{u}_k)\cdot\partial_t(\sqrt{\alpha_k}\ve{u}_k)\, \mathrm{d}\Omega\, ,\\
  &B_k = \rho_k\int_{\Omega}\ve{u}_k\cdot\bigg[\alpha_k\ve{u}_k\cdot\nabla\ve{u}_k  + \frac{1}{2} \nabla\cdot(\alpha_k\ve{u}_k)\ve{u}_k\bigg]\, \mathrm{d}\Omega\, \\
  &C_k = -\int_{\Omega}\ve{u}_k\cdot\left[\nabla\cdot(2\alpha_k\mu_k\nabla^{\mathrm{s}}\ve{u}_k)\right] \mathrm{d}\Omega\, , \\
  &D_k = \int_{\Omega}\ve{u}_k\cdot\sum_{l=1}^{M}\gamma_{kl}(\ve{u}_k-\ve{u}_l)\, \mathrm{d}\Omega\, , \\
  &E_k = \int_{\Omega}\alpha_k\ve{u}_k\cdot\nabla p \, \mathrm{d}\Omega\, .
\end{align*}

The first term can be rewritten as
\begin{align*}
  A_k &= \frac{1}{2}\rho_k\int_{\Omega}\partial_t(|\sqrt{\alpha_k}\ve{u}_k|^2)\, \mathrm{d}\Omega\\
  &= \frac{\mathrm{d}}{\mathrm{d} t}\left(\frac{1}{2}\rho_k\|\sqrt{\alpha_k}\ve{u}_k\|^2\right),
\end{align*}
with $| \cdot |$ denoting the Euclidean norm. For the convective term, we have
\begin{align*} \ve{u}_k\cdot\bigg[\alpha_k\ve{u}_k\cdot\nabla\ve{u}_k  + \frac{1}{2} \nabla\cdot(\alpha_k\ve{u}_k)\ve{u}_k\bigg] &= \alpha_k\ve{u}_k\cdot[\ve{u}_k\cdot\nabla\ve{u}_k]  + \nabla\cdot(\alpha_k\ve{u}_k)\frac{|\ve{u}_k|^2}{2}\\
&=\alpha_k\ve{u}_k\cdot\nabla\bigg(\frac{|\ve{u}_k|^2}{2}\bigg)  + (\ve{u}_k\cdot\nabla\alpha_k+\alpha_k\nabla\cdot\ve{u}_k)\frac{|\ve{u}_k|^2}{2}\\
&=\frac{1}{2}\ve{u}_k\cdot\left[\alpha_k\nabla\big(|\ve{u}_k|^2\big)+|\ve{u}_k|^2\nabla\alpha_k\right]  + \frac{\alpha_k|\ve{u}_k|^2}{2}\nabla\cdot\ve{u}_k\\
&=\frac{1}{2}\left[\ve{u}_k\cdot\nabla\big(\alpha_k|\ve{u}_k|^2\big)  + |\ve{u}_k|^2\nabla\cdot\ve{u}_k\right]\\
&= \frac{1}{2}\nabla\cdot(\alpha_k|\ve{u}_k|^2\ve{u}_k)\, ,
\end{align*}
so that 
\begin{align*}
    B_k &= \frac{\rho_k}{2}\int_{\Omega}\nabla\cdot(\alpha_k|\ve{u}_k|^2\ve{u}_k) \, \mathrm{d}\Omega \\
    &=\frac{\rho_k}{2}\int_{\partial\Omega}\ve{n}\cdot(\alpha_k|\ve{u}_k|^2\ve{u}_k) \, \mathrm{d}\Gamma\\
    &= 0\, ,
\end{align*}
due to the no-slip boundary condition. For the viscous term, integration by parts yields 
\begin{align*}
    C_k &= \int_{\Omega}2\alpha_k\mu_k\nabla^{\mathrm{s}}\ve{u}_k:\nabla^{\mathrm{s}}\ve{u}_k\, \mathrm{d}\Omega - \int_{\partial\Omega} \ve{u}_k\cdot\left[(2\alpha_k\mu_k\nabla^{\mathrm{s}}\ve{u}_k)\ve{n}\right] \mathrm{d}\Gamma\\
    &= 2\|\sqrt{\mu_k\alpha_k}\,\nabla^{\mathrm{s}}\ve{u}_k\|^2\, .
\end{align*}

To compute the total dissipation caused by the drag forces over all phases, we write
\begin{align*}
    \sum_{k=1}^M \left[\ve{u}_k\cdot\sum_{l=1}^{M}\gamma_{kl}(\ve{u}_k-\ve{u}_l)\right] &= \sum_{k=1}^M\sum_{l=1}^{M} \ve{u}_k\cdot\gamma_{kl}(\ve{u}_k-\ve{u}_l)\\
    &=\frac{1}{2}\sum_{k=1}^M\sum_{l=1}^{M} \gamma_{kl}\ve{u}_k\cdot(\ve{u}_k-\ve{u}_l) + \frac{1}{2}\sum_{k=1}^M\sum_{l=1}^{M} \gamma_{kl}\ve{u}_k\cdot(\ve{u}_k-\ve{u}_l)\\
    &=\frac{1}{2}\sum_{k=1}^M\sum_{l=1}^{M} \gamma_{kl}\ve{u}_k\cdot(\ve{u}_k-\ve{u}_l) + \frac{1}{2}\sum_{k=1}^M\sum_{l=1}^{M} \gamma_{lk}\ve{u}_l\cdot(\ve{u}_l-\ve{u}_k)\\
    &=\frac{1}{2}\sum_{k=1}^M\sum_{l=1}^{M}\gamma_{kl}\left[ \ve{u}_k\cdot(\ve{u}_k-\ve{u}_l) -  \ve{u}_l\cdot(\ve{u}_k-\ve{u}_l)\right]\\
    &=\frac{1}{2}\sum_{k=1}^M\sum_{l=1}^{M}\gamma_{kl}|\ve{u}_k-\ve{u}_l|^2\, ,
\end{align*}
where we have used that $\gamma_{kl} = \gamma_{lk}$. Hence: 
\begin{align*}
    \sum_{k=1}^M D_k = \frac{1}{2}\sum_{k=1}^M\sum_{l=1}^{M}\big\|\sqrt{\gamma_{kl}}(\ve{u}_k-\ve{u}_l)\big\|^2\, .
\end{align*}

For the pressure terms $E_k$, summing over all phases and integrating by parts yields 
\begin{align*}
   \sum_{k=1}^M E_k &= \int_{\Omega}\nabla p\cdot\sum_{k=1}^M(\alpha_k\ve{u}_k) \, \mathrm{d}\Omega\\
   &= -\int_{\Omega}p\nabla\cdot\sum_{k=1}^M(\alpha_k\ve{u}_k) \, \mathrm{d}\Omega + \int_{\partial\Omega}p\ve{n}\cdot\sum_{k=1}^M(\alpha_k\ve{u}_k) \, \mathrm{d}\Gamma\\
   &= -\int_{\Omega}p\nabla\cdot\bar{\ve{u}} \, \mathrm{d}\Omega\\
   &= 0\, ,
\end{align*}
where we have used once again the zero boundary conditions. 

Finally, collecting the terms for all phases gives us
\begin{align*}
   \sum_{k=1}^{M}\left[\frac{\mathrm{d}}{\mathrm{d} t}\left(\frac{1}{2}\rho_k\|\sqrt{\alpha_k}\ve{u}_k\|^2\right) + 2\|\sqrt{\mu_k\alpha_k}\,\nabla^{\mathrm{s}}\ve{u}_k\|^2 + \frac{1}{2}\sum_{l=1}^{M}\big\|\sqrt{\gamma_{kl}}(\ve{u}_k-\ve{u}_l)\big\|^2 \right] = 0\, ,
\end{align*}
which, when integrated over time, yields the conservation equation 
\begin{equation}
\begin{split}
    &\sum_{k=1}^{M}\left\lbrace\frac{1}{2}\rho_k\|\sqrt{\alpha_k}\ve{u}_k\|^2(t) + \int_{0}^t\left[2\|\sqrt{\mu_k\alpha_k}\,\nabla^{\mathrm{s}}\ve{u}_k\|^2 + \frac{1}{2}\sum_{l=1}^{M}\big\|\sqrt{\gamma_{kl}}(\ve{u}_k-\ve{u}_l)\big\|^2 \right]\mathrm{d}t\right\rbrace\\
    &= - \sum_{k=1}^{M}\frac{1}{2}\rho_k\|\sqrt{\alpha_k}\ve{u}_k\|^2(0)\, .\label{energyConservation}
    \end{split}
\end{equation}
This shows that, for a confined system (no volume or boundary energy sources), the decrease in kinetic energy over time is equal to the work done by the viscous and drag forces. While some fully implicit monolithic schemes may inherit this property at the discrete level---which leads to numerical stability---constructing a stable fractional-step method whose energy inequality resembles \eqref{energyConservation} is much more challenging. The design and analysis of such a method is the main contribution of this article. 

\subsection{Mass}
The volume fraction equations are a delicate matter in terms of conservation and stability, since the individual velocities $\ve{u}_k$ are not divergence-free. For a confined flow, integrating Eq.~\eqref{alpha1} over space and time yields
\begin{align*}
    0 &= \int_0^t\int_{\Omega}\left[\partial_t\alpha_k + \nabla\cdot(\alpha_k\ve{u}_k)\,\right] \mathrm{d}\Omega\,\mathrm{d}t\\
    &= \int_0^t\frac{\mathrm{d}}{\mathrm{d}t}\left(\int_{\Omega}\alpha_k\,  \mathrm{d}\Omega\right)\mathrm{d}t + \int_0^t\left[\int_{\partial\Omega}\ve{n}\cdot(\alpha_k\ve{u}_k)\,\mathrm{d}\Gamma\right]\mathrm{d}t\\
    &= \int_{\Omega}\alpha_k(t)\,  \mathrm{d}\Omega - \int_{\Omega}\alpha_k(0)\,  \mathrm{d}\Omega\, ,
\end{align*}
which, when multiplied by the constant density $\rho_k$, shows conservation of mass for each phase over time. This property, however, does not translate to (numerical) stability, since the integral of $\alpha_k$ is not a norm. Traditional energy-like estimates, e.g.~using $L^2$ norms, require the velocity field to be solenoidal, which is not the case here. The numerical implications of this fact will be discussed in Section \ref{sec_convection}.

\section{Temporal discretisation}\label{sec_schemes}
We shall denote approximate or discrete values of quantities at different time steps with superscripts: $\ve{u}^{n}_k$, for example, denotes the approximation of $\ve{u}_k$ at $t=t_n = n\tau$, with the time-step size $\tau$ assumed as constant, for simplicity. To construct the time-stepping scheme, we will approximate the temporal derivatives with first-order finite differences, while using first-order extrapolations (e.g., $\ve{u}_k^{n+1}\approx\ve{u}_k^{n}$) to avoid nonlinearities and terms that would induce undesired coupling if treated implicitly.

For the pressure-momentum system, we will present a fractional-step method to compute each velocity $\ve{u}_k$ and the pressure $p$ all separately. As a crucial feature for robustness, our method can be derived and proved stable \textsl{independently} of how the volume fractions are updated---provided that they, too, are evolved in a stable, positivity-preserving way. Thus, for now, we shall assume a generic volume fraction scheme $f$ that computes $\alpha_k^{n+1}$ based on $\alpha_k^n$ and on the (extrapolated) velocity $\ve{u}_k^n$:
\begin{equation*}
    \alpha_k^{n+1} = f(\tau,\alpha_k^{n},\ve{u}_k^n)\, .
\end{equation*}
Different formulations will be discussed later on in Section \ref{sec_convection}.

\subsection{Convective terms}
For the inertial terms, we first write
\begin{align}
&\left.\left\lbrace\sqrt{\alpha_k}\,\partial_t(\sqrt{\alpha_k}\ve{u}_k) + \alpha_k\ve{u}_k\cdot\nabla\ve{u}_k + \frac{1}{2}[\nabla\cdot(\alpha_k\ve{u}_k)]\ve{u}_k\right\rbrace\right|_{t=t_{n+1}} \nonumber\\
			&\approx \frac{\sqrt{\alpha_k^{n+1}}}{\tau}\left(\sqrt{\alpha_k^{n+1}}\ve{u}_k^{n+1}-\sqrt{\vphantom{\alpha_k^{n+1}}\alpha_k^{n}}\,\ve{u}_k^n\right) + \alpha_k^{n+1}\ve{u}_k^n\cdot\nabla\ve{u}_k^{n+1}  + \frac{1}{2}\left[ \nabla\cdot(\alpha_k^{n+1}\ve{u}_k^n)\right]\ve{u}_k^{n+1}\nonumber\\
			&=  \frac{1}{\tau}\left(\alpha_k^{n+1}\ve{u}_k^{n+1}-\sqrt{\alpha_k^{n+1}\alpha_k^{n}}\,\ve{u}_k^n\right) + \alpha_k^{n+1}\ve{u}_k^n\cdot\nabla\ve{u}_k^{n+1}  + \frac{1}{2}\left[ \nabla\cdot(\alpha_k^{n+1}\ve{u}_k^n)\right]\ve{u}_k^{n+1}\, .\label{convectiveBDF1}
\end{align}
Although one may alternatively treat the convective term implicitly, our choice here is the IMEX approach \eqref{convectiveBDF1}. The main reason is to circumvent the nonlinearity while still maintaining unconditional stability, as proved in Section \ref{sec_analysis}; in contrast, fully explicit treatments are known to introduce CFL conditions, even for the single-phase case \cite{Burman2023,Burman2024}. Now, for the discretised acceleration term on the right-hand side of \eqref{convectiveBDF1}, we introduce the first-order extrapolation
    \begin{align}
        \alpha_k^{n+1}\ve{u}_k^{n+1}-\sqrt{\alpha_k^{n+1}\alpha_k^{n}}\,\ve{u}_k^n &= \frac{\alpha_k^{n+1}+\alpha_k^{n+1}}{2}\ve{u}_k^{n+1}-\sqrt{\alpha_k^{n+1}\alpha_k^{n}}\,\ve{u}_k^n\nonumber \\
        &\approx \frac{\alpha_k^{n+1}+\alpha_k^{n}}{2}\ve{u}_k^{n+1}-\sqrt{\alpha_k^{n}\alpha_k^{n}}\,\ve{u}_k^n\nonumber\\
        &= \frac{\alpha_k^{n+1}+\alpha_k^{n}}{2}\ve{u}_k^{n+1}-\alpha_k^{n}\,\ve{u}_k^n\, .\label{Salgado}
    \end{align}
    This reformulation may seem somewhat unnatural, yet it has a (rather technical) numerical purpose: in short, the original term in \eqref{convectiveBDF1}, when tested against $\ve{u}_{k}^{n+1}$, would give us
\begin{align*}    &\left\langle\alpha_k^{n+1}\ve{u}_k^{n+1}-\sqrt{\alpha_k^{n+1}\alpha_k^{n}}\,\ve{u}_k^n,\ve{u}_{k}^{n+1} \right\rangle \\
    &= \frac{1}{2}\left(\left\|\sqrt{\alpha_k^{n+1}}\,\ve{u}_k^{n+1}\right\|^2 - \left\|\sqrt{\alpha_k^{n}}\,\ve{u}_k^{n}\right\|^2 + \left\|\sqrt{\alpha_k^{n+1}}\,\ve{u}_k^{n+1} - \sqrt{\alpha_k^{n}}\,\ve{u}_k^{n}\right\|^2\right),
\end{align*}
 to \eqref{identityBDF1}, whereas our reformulated version yields (see \eqref{estAcc} later on)
\begin{align*}
&\left\langle\frac{\alpha_k^{n+1}+\alpha_k^{n}}{2}\ve{u}_k^{n+1}-\alpha_k^{n}\,\ve{u}_k^n,\ve{u}_{k}^{n+1} \right\rangle \\
    &= \frac{1}{2}\left(\left\|\sqrt{\alpha_k^{n+1}}\,\ve{u}_k^{n+1}\right\|^2 - \left\|\sqrt{\alpha_k^{n}}\,\ve{u}_k^{n}\right\|^2 + \left\|\sqrt{\alpha_k^{n}}(\ve{u}_k^{n+1} - \ve{u}_k^{n})\right\|^2\right).
\end{align*}
In comparison to the term $\sqrt{\alpha_k^{n+1}}\,\ve{u}_k^{n+1} - \sqrt{\vphantom{\alpha_k^{n+1}}\alpha_k^{n}}\,\ve{u}_k^{n}$, the modified one $\sqrt{\vphantom{\frac{\cdot}{\cdot}}\alpha_k^{n}}(\ve{u}_k^{n+1} - \ve{u}_k^{n})$ is more convenient for controlling the explicit drag, since both velocities have the same factor ($\sqrt{\vphantom{\frac{\cdot}{\cdot}}\alpha_k^{n}}$), similarly as we have for the relative velocities $\ve{u}_k-\ve{u}_l$ appearing in the drag forces. This will be important to guarantee temporal stability for the explicit drag and will hopefully become clearer once we start the analysis.

    \subsection{Viscous terms}
        Although the viscous term $2\alpha_k\mu_k\nabla^{\mathrm{s}}\ve{u}_k$ itself does not couple the $M$ phase velocities, if made fully implicit it has the (admittedly minor) disadvantage of coupling the $d$ spatial components of $\ve{u}_k$, due to the transpose gradient in $\nabla^{\mathrm{s}}\ve{u}_k=\nicefrac{1}{2}(\nabla\ve{u}_k + \nabla^{\top}\ve{u}_k)$. To enable a component-wise solution of $\ve{u}_k$ by solving scalar subproblems, we propose the first-order extrapolation 
        \begin{align*}
(2\alpha_k\mu_k\nabla^{\mathrm{s}}\ve{u}_k)|_{t=t_{n+1}}&= \mu_k\alpha_k^{n+1}(\nabla\ve{u}_k^{n+1}+\nabla^{\top}\ve{u}_k^{n+1})\\
&= \mu_k\alpha_k^{n+1}\nabla\ve{u}_k^{n+1}+\mu_k\sqrt{\alpha_k^{n+1}\alpha_k^{n+1}}\,\nabla^{\top}\ve{u}_k^{n+1}\\
&\approx \mu_k\alpha_k^{n+1}\nabla\ve{u}_k^{n+1}+\mu_k\sqrt{\alpha_k^{n+1}\alpha_k^{n}}\,\nabla^{\top}\ve{u}_k^{n}\, .
        \end{align*}
We will show later that this IMEX treatment, first proposed in the context of single-phase variable-viscosity flows \cite{Barrenechea2024}, is stable in time. 

\begin{remark}
There are applications where some of the viscosities $\nu_k$ may depend on $\ve{u}_k$ and/or on $p$, as in dense granular flows \cite{Gesenhues2021}. In that scenario, to avoid viscous nonlinearities, one can use a ``nearly implicit'' approach extrapolating only the viscosities:
\begin{align*}
&\left.( 2\rho\nu_k\alpha_k\nabla^{\mathrm{s}}\ve{u}_k)\right|_{t=t_{n+1}} \approx 2\rho_k\nu_k^n\alpha_k^{n+1}\nabla^{\mathrm{s}}\ve{u}_k^{n+1}\, .
\end{align*}
Making the viscosity explicit as above, which is compatible with the $\mathcal{O}(\tau)$ consistency of our framework, yields essentially the same stability as a fully implicit treatment.
\end{remark}
	
\subsection{Numerical method}
To fully decouple the phases, we must treat both the pressure and the drag forces explicitly. For the pressure, we propose a fractional-step scheme that can be seen as a multiphase extension of a recently proposed projection method for diffuse-interface flows \cite{Deteix2022}. A key idea here is to use a first-order extrapolation of the pressure term:
\begin{align*}
    \alpha_k^{n+1}\nabla p^{n+1} 
    \approx \sqrt{\alpha_k^{n+1}\alpha_k^{n}}\,\nabla p^{n}\, .
\end{align*}
We also introduce a modification of the drag coefficients:
\begin{align*}
   \tilde{\gamma}_{kl}^{n} := \min \left\lbrace\gamma_{kl}^n,D\right\rbrace,
\end{align*}
where $D>0$ is a large, but finite constant. This is a technical assumption needed for the analysis, with probably little practical relevance.  
		
We will consider Dirichlet boundary conditions for the momentum equations \eqref{momentum2}:
	\begin{align*}
		\ve{u}_k = \mathbf{u}_k^D \ \ \text{on} \ \partial\Omega\times(0,T]&\, .
	\end{align*}
The scheme presented below can also be used for the free-slip case ($\ve{n}\cdot\ve{u}_k = 0$, with the tangential component left ``free''). For open (Neumann) boundaries, however, a modified scheme would be required, which is a challenging task that will not be addressed herein. For the volume fractions, boundary conditions are only needed (and allowed) for $\alpha_k$ whenever and wherever there is boundary inflow ($\ve{n}\cdot\ve{u}_k < 0$).

Finally, we propose the following \textsl{pressure-correction} method:
	\begin{itemize}
		\item \textbf{Step 0:} Set $\hat{\mathbf{u}}_k^0 = \ve{u}_k^0$ and provide (an approximation for) $p^0$ (see Section \ref{sec_p0}).
		
		\item \textbf{Step 1:} Update each volume fraction $\alpha_k$ using the scheme of choice.
		
		\item \textbf{Step 2:} Find $\ve{u}^{n+1}_k$ by solving the linear advection-diffusion-reaction problem
		\begin{flalign}
			\begin{cases}
	\rho_k\Big[\frac{\alpha_k^{n+1}+\alpha_k^{n}}{2\tau}\ve{u}_k^{n+1}   + \alpha_k^{n+1}\ve{u}_k^n\cdot\nabla\ve{u}_k^{n+1} + \frac{1}{2}\nabla\cdot(\alpha_k^{n+1}\ve{u}_k^n)\ve{u}_k^{n+1}\Big]  -\nabla\cdot(\mu_k\alpha^{n+1}_k\nabla\ve{u}^{n+1}_k)  \\
			=   \rho_k\frac{\alpha_k^{n}}{\tau}\hat{\mathbf{u}}_k^{n}  +  \nabla\cdot\Big(\mu_k\sqrt{\alpha^{n+1}_k\alpha^{n}_k}\,\nabla^{\top}\ve{u}^{n}_k\Big)  -\sqrt{\alpha_k^{n+1}\alpha_k^{n}}\,\nabla p^{n} \\
             + \rho_k\alpha_k^{n+1}\ve{g}_k^{n+1} - \sum_{l=1}^{M}\tilde{\gamma}_{kl}^{n}(\hat{\mathbf{u}}_k^{n}-\hat{\mathbf{u}}_l^{n})
            \\
				\ve{u}_k^{n+1}|_{\partial\Omega} = \left.\mathbf{u}_k^{D}\right|_{t=t_{n+1}}
			\end{cases} && \label{FS-momentum}
		\end{flalign}
		for each $k=1,\dots,M$ (the viscous term could alternatively be treated implicitly). Notice that the so-called \textsl{end-of-step velocities} $(\hat{\mathbf{u}}^n_k,\hat{\mathbf{u}}^n_l)$ are used in the first and last terms on the right-hand side of the momentum equation. The reason for using these projected velocities to evaluate the drag will become clear during the stability analysis, where it is shown that the fractional-step scheme provides $L^2$ control over $\hat{\mathbf{u}}^n_k$ (as usual in projection methods \cite{Deteix2022,Guermond2000}) instead of $\ve{u}^n_k$.
		
		\item \textbf{Step 3:} Find the pressure $p^{n+1}$ through the Poisson problem
		\begin{flalign}
			\begin{cases}
\nabla\cdot\left[\Big(\sum_{k=1}^{M}\rho_k^{-1}\alpha_k^{n+1}\Big)\nabla p^{n+1}\right] 
\\= \nabla\cdot\left[\Big(\sum_{k=1}^{M}\rho_k^{-1}\sqrt{\alpha_k^{n+1}\alpha_k^{n}}\Big)\nabla p^{n}\right] +\frac{1}{\tau}\nabla\cdot\sum_{k=1}^{M}\alpha_k^{n+1}\ve{u}_k^{n+1}\\
\left.\Big(\sum_{k=1}^{M}\rho_k^{-1}\alpha_{k}^{n+1}\Big) \partial_{\vesmall{n}}p^{n+1}\right|_{\partial\Omega} = \Big(\sum_{k=1}^{M}\rho_k^{-1}\sqrt{\alpha_{k}^{n+1}\alpha_{k}^{n}}\left.\Big)\partial_{\vesmall{n}}p^{n}\right|_{\partial\Omega}
\end{cases}&& \label{PPE}
\end{flalign}
\item \textbf{Step 4:} Update each end-of-step velocity via 
\begin{flalign}
			\hat{\mathbf{u}}_k^{n+1} &= \ve{u}_k^{n+1} + \frac{\tau}{\rho_k}\left(\sqrt{\frac{\alpha_k^n}{\alpha_k^{n+1}}}\,\nabla p^{n}-\nabla p^{n+1}\right) \, . && \label{projectedU}
\end{flalign}
\end{itemize}
	
    Notice that combining Eqs.~\eqref{PPE} and \eqref{projectedU} yields
		\begin{align}
\nabla\cdot\sum_{k=1}^{M}\alpha_k^{n+1}\hat{\mathbf{u}}_k^{n+1} = 0 \ \ \text{and}
			\ \ \ve{n}\cdot\sum_{k=1}^{M}\alpha_k^{n+1}(\ve{u}_k^{n+1}-\hat{\mathbf{u}}_k^{n+1})|_{\partial\Omega} = 0\, ,\label{divEndOfStep}
		\end{align}
that is, the volume-averaged end-of-step velocity is indeed a divergence-free projection. This, along with the derivation of the fractional-step scheme, is shown in \ref{sec_appendix}.
    
Not only are the steps for each phase completely independent and decoupled, but they are also linear. After updating each volume fraction, we solve $M$ advection-diffusion-reaction problems to compute all the $\ve{u}_k^{n+1}$, then a single pressure Poisson equation (PPE). Moreover, due to the IMEX treatment of the convective and viscous terms, each of the momentum equations can be broken into $d$ (2 or 3) \textsl{scalar} advection-diffusion-reaction problems (one for each spatial component of $\ve{u}_k^{n+1}$). These algorithmic properties enable a simple implementation, since each subproblem is a (variable-coefficient) linear transport equation. Another often understated fact is the additional robustness brought by the linearity of the scheme. Here, there is no issue of fixed-point convergence: for any $\tau$ the solution of \textsl{every} time step converges---provided, of course, that the discretisation scheme itself is stable, which we will soon prove.


\begin{remark}
From Eq.~\eqref{PPE} we have that
			\begin{align*}		\partial_{\vesmall{n}}p^{n+1} &= \left(\frac{\sum_{k=1}^{M}\rho_k^{-1}\sqrt{\alpha_{k}^{n+1}\alpha_{k}^{n}}}{\sum_{k=1}^{M}\rho_k^{-1}\alpha_{k}^{n+1}}\right)\partial_{\vesmall{n}} p^{n}\\
            &=\dots =\left(\prod_{m=1}^{n+1}\frac{\sum_{k=1}^{M}\rho_k^{-1}\sqrt{\alpha_{k}^{m+1}\alpha_{k}^{m}}}{\sum_{k=1}^{M}\rho_k^{-1}\alpha_{k}^{m+1}}\right)\partial_{\vesmall{n}} p^{0}
				\ \ \text{on} \ \partial\Omega\, ,
			\end{align*}
which means that, as most standard pressure-correction schemes, ours also induces an artificial pressure boundary condition. This tends to create a numerical boundary layer for large time-step sizes. Similarly, the end-of-step velocities $\hat{\mathbf{u}}_k$ do not exactly retain the Dirichlet boundary conditions enforced on $\ve{u}_k$. These artifacts prevent standard projection methods from being high-order accurate.
\end{remark}

\subsection{Interpretation as a penalty-like method}\label{sec_penalty}
While convergence analysis is not in the present scope, we will sketch how our projection scheme can be reinterpreted as a penalty-like method, which will also offer insights into its order of consistency. Evaluating Eq.~\eqref{projectedU} at $t=t_n$ allows us to eliminate the term $\rho_k\alpha_k^n\hat{\mathbf{u}}_k^n$ from Eq.~\eqref{FS-momentum}, leaving us with
\begin{align*}
    &\rho_k\Big[\frac{\alpha_k^{n+1}+\alpha_k^{n}}{2\tau}\ve{u}_k^{n+1}   + \alpha_k^{n+1}\ve{u}_k^n\cdot\nabla\ve{u}_k^{n+1} + \frac{1}{2}\nabla\cdot(\alpha_k^{n+1}\ve{u}_k^n)\ve{u}_k^{n+1}\Big]  -\nabla\cdot(\mu_k\alpha^{n+1}_k\nabla\ve{u}^{n+1}_k)  \\
			&=   \rho_k\frac{\alpha_k^{n}}{\tau}\ve{u}_k^{n}  +  \nabla\cdot\Big(\mu_k\sqrt{\alpha^{n+1}_k\alpha^{n}_k}\,\nabla^{\top}\ve{u}^{n}_k\Big) + \rho_k\alpha_k^{n+1}\ve{g}_k^{n+1} + \sum_{l=1}^{M}\tilde{\gamma}_{kl}^{n}(\hat{\mathbf{u}}_l^{n}-\hat{\mathbf{u}}_k^{n})\\
             &-\left(\sqrt{\alpha_k^{n+1}\alpha_k^{n}}\,\nabla p^{n} + \alpha_k^n\nabla p^n -\sqrt{\alpha_k^{n}\alpha_k^{n-1}}\,\nabla p^{n-1} 
             \right) .
\end{align*}
Although maybe not obvious at first sight, the pressure terms on the right-hand side form a second-order extrapolation of $-\alpha_k^{n+1}\nabla p^{n+1}$. To show that, let us denote $\sqrt{\alpha_k}=\varphi$ and $\sqrt{\alpha_k}\,\nabla p = \ve{a}$, so that a Taylor expansion about $t_n$ gives us
\begin{align*}
&\sqrt{\alpha_k^{n+1}\alpha_k^{n}}\,\nabla p^{n}  + \alpha_k^n\nabla p^n
              -\sqrt{\alpha_k^{n}\alpha_k^{n-1}}\,\nabla p^{n-1} \\
             &= \varphi^{n+1}\ve{a}^{n} + \varphi^{n}\ve{a}^{n} - \varphi^{n}\ve{a}^{n-1}\\
             &= \big[\varphi^{n}+\tau(\partial_t\varphi)|_{t=t_{n}} + \mathcal{O}(\tau^2)\big]\ve{a}^{n} + \varphi^{n}\ve{a}^{n} - \varphi^{n}\big[\ve{a}^{n} - \tau(\partial_t\ve{a})|_{t=t_{n}} + \mathcal{O}(\tau^2)\big]\\
             &= \varphi^{n}\ve{a}^{n}+\tau\big[(\ve{a}\partial_t\varphi)|_{t=t_{n}}  + (\varphi\partial_t\ve{a})|_{t=t_{n}}\big]  + \mathcal{O}(\tau^2)\\
             &= (\varphi\ve{a})|_{t=t_n}+\tau[\partial_t(\varphi\ve{a})]|_{t=t_n} + \mathcal{O}(\tau^2)\\
             &= (\varphi\ve{a})|_{t=t_{n+1}}+ \mathcal{O}(\tau^2)\\
             &= \alpha_k^{n+1}\nabla p^{n+1} + \mathcal{O}(\tau^2)\, .
\end{align*}
Moreover, denoting $\bar{\ve{u}}^{n+1} = \sum_{k=1}^{M}\alpha_k^{n+1}\ve{u}_k^{n+1}$, the PPE \eqref{PPE} can be rewritten as 
\begin{align*}
    \nabla\cdot\bar{\ve{u}}^{n+1} &= \tau\nabla\cdot\sum_{k=1}^M\rho_k^{-1}\sqrt{\alpha_k^{n+1}}\left(\sqrt{\alpha_k^{n+1}}\,\nabla p^{n+1}-\sqrt{\vphantom{\alpha_k^{n+1}}\alpha_k^{n}}\,\nabla p^{n}\right)\\
    &= \tau^2\nabla\cdot\sum_{k=1}^M\rho_k^{-1}\sqrt{\alpha_k^{n+1}}\left(\frac{\sqrt{\alpha_k^{n+1}}\,\nabla p^{n+1}-\sqrt{\vphantom{\alpha_k^{n+1}}\alpha_k^{n}}\,\nabla p^{n}}{\tau}\right)\\
    &= \tau^2\nabla\cdot\sum_{k=1}^M\rho_k^{-1}\big[\sqrt{\alpha_k}\,\partial_t(\sqrt{\alpha_k}\,\nabla p)\big]|_{t=t_{n+1}} + \mathcal{O}(\tau^3)\, ,
\end{align*}
that is, as an $\mathcal{O}(\tau^2)$ perturbation (the term proportional to $\tau^2$ dominates the cubic one) of the incompressibility constraint on the mean velocity. In other words, the fractional-step method itself---which need not necessarily be treated in the IMEX fashion we have proposed---can be seen as a (first-order) discretisation of the $\mathcal{O}(\tau^2)$-perturbed system 
\begin{flalign*}
    \begin{cases}
\text{Momentum eqs.~\eqref{momentum2}}\, , \\
\nabla\cdot\bar{\ve{u}} -\tau^2\sum_{k=1}^{M} \rho_k^{-1}\nabla\cdot\left[\sqrt{\alpha_k}\,\partial_t(\sqrt{\alpha_k}\,\nabla p)\right] = 0\, .
    \end{cases} &&
\end{flalign*}
This indicates that it should be possible to extend our method to order two in time---although guaranteeing stability of such a scheme could be considerably more difficult.

\begin{remark}
  In the single-phase case, the perturbation to the incompressibility constraint would reduce to the well-known term $-\tau^2\partial_t\Delta p$, as in the incremental pressure-correction method \cite{Guermond2006}. 
\end{remark}

\section{Weak formulation}
	The method presented in Steps 0--4 above is meant as a general time-stepping scheme that, in principle, may be combined with any spatial discretisation of choice. For the case of finite element methods, each subproblem needs to be written in weak form. Without dwelling with finite element formalism, let us simply assume $X_h$ and $Y_h$ as globally continuous, $H^1$-conforming finite element spaces for velocities and pressure, respectively. Dirichlet or inlet conditions will be omitted, for concision, so it is implied that test and trial functions satisfy the respective boundary conditions.

  \subsection{Basic steps}  
    The weak subproblems are:
	\begin{itemize}
		\item \textbf{Step 1:} Compute each volume fraction $\alpha_k^{n+1}$ using the method of choice.
		
		\item \textbf{Step 2:} for each $k=1,\dots,M$, find $\ve{u}_k^{n+1}\in X_h$ such that 
\begin{align}
				&\rho_k\left\langle \frac{\alpha_k^{n+1}+\alpha_k^{n}}{2\tau}\ve{u}_k^{n+1}   + \alpha_k^{n+1}\ve{u}_k^{n}\cdot\nabla\ve{u}_k^{n+1} + \frac{1}{2}[\nabla\cdot(\alpha_k^{n+1}\ve{u}_k^n)]\ve{u}_k^{n+1},\ve{v} \right\rangle+ \langle \mu_k\alpha_k^{n+1}\nabla\ve{u}_k^{n+1},\nabla\ve{v}\rangle \nonumber \\ 
				&=  \frac{\rho_k}{\tau}\langle\alpha_k^{n}\hat{\mathbf{u}}_k^{n},\ve{v}\rangle + \left\langle p^{n},\nabla\cdot\left(\sqrt{\alpha_k^{n+1}\alpha_k^{n}}\,\ve{v}\right)\right\rangle - \left\langle \mu_k\sqrt{\alpha_k^{n+1}\alpha_k^{n}}\,\nabla^{\top}\ve{u}_k^{n},\nabla\ve{v}\right\rangle \label{weakMomentum}\\
                &+ \rho_k\langle\alpha_k^{n+1}\ve{g}_k^{n+1},\ve{v}\rangle - \sum_{l=1}^{M}\left\langle\tilde{\gamma}_{kl}^{n}(\hat{\mathbf{u}}_k^{n}-\hat{\mathbf{u}}_l^{n}),\ve{v}\right\rangle \nonumber    
\end{align} 
		for all $\ve{v}\in X_h$.
		
		\item \textbf{Step 3:} find $p^{n+1}\in Y_h$ such that 
		\begin{equation}
			\left\langle\nabla p^{n+1},\left(\sum_{k=1}^{M}\frac{\alpha_k^{n+1}}{\rho_k}\right) \nabla q \right\rangle = \left\langle\nabla p^{n},\left(\sum_{k=1}^{M}\frac{\sqrt{\alpha_k^{n+1}\alpha_k^{n}}}{\rho_k}\right)\nabla q \right\rangle -\frac{1}{\tau}\sum_{k=1}^{M}\left\langle\nabla\cdot(\alpha_k^{n+1}\ve{u}_k^{n+1}),q\right\rangle\label{weakPPE}
		\end{equation}
		for all $q\in Y_h$.
		
		\item \textbf{Step 4:} for each $k=1,\dots,M$, find $\hat{\mathbf{u}}_k^{n+1}\in X_h$ such that 
		\begin{equation}
\left\langle\alpha_k^{n+1}\hat{\mathbf{u}}_k^{n+1},\mathbf{v}\right\rangle = \left\langle\alpha_k^{n+1}\ve{u}_k^{n+1},\mathbf{v}\right\rangle + \frac{\tau}{\rho_k}\left\langle\sqrt{\vphantom{\alpha_k^{n+1}}\alpha_k^{n}\,}\,\nabla p^{n}-\sqrt{\alpha_k^{n+1}}\,\nabla p^{n+1},\sqrt{\alpha_k^{n+1}}\,\mathbf{v}\right\rangle \label{weakUhat}
		\end{equation}
		for all $\mathbf{v}\in X_h$.
	\end{itemize}

\subsection{Initialising the pressure}\label{sec_p0}
	For incompressible flows, the initial pressure $p^0$ is not part of the required problem data. However, explicit or IMEX schemes typically require $p^0$ for initialisation. The approach we describe next assumes non-penetrating boundaries, $\ve{n}\cdot\ve{u}_k|_{\partial\Omega}=0$, but more general setups can also be tackled with minor modifications. We start by dividing each momentum equation by \eqref{momentum} by $\rho_k$ and adding them all together:
	\begin{align}
\left(\sum_{k=1}^{M}\frac{\alpha_k}{\rho_k}\right)\nabla p+ \sum_{k=1}^{M}\left[\partial_t(\alpha_k\ve{u}_k) - \ve{f}_k\right] = \ve{0}  \, ,\label{gradientEq}
	\end{align}
	where $\ve{f}_k := -\nabla\cdot(\alpha_k\ve{u}_k\otimes\ve{u}_k) + \nabla\cdot(2\nu_k\alpha_k\nabla^{\mathrm{s}}\ve{u}_k) + \alpha_k\ve{g}_k - \rho_k^{-1}\sum_{l=1}^{M}\gamma_{kj}(\ve{u}_k-\ve{u}_l)$. We then obtain a Poisson equation for the pressure by taking the divergence of Eq.~\eqref{gradientEq}:
	\begin{align*}
		\nabla\cdot\left[\left(\sum_{k=1}^{M}\frac{\alpha_k}{\rho_k}\right)\nabla p\right] - \nabla\cdot\sum_{k=1}^{M}\ve{f}_k &= -\nabla\cdot\sum_{k=1}^{M}\partial_t(\alpha_k\ve{u}_k)\\
		&= -\partial_t\left(\nabla\cdot\sum_{k=1}^{M}\alpha_k\ve{u}_k\right)\\
		&= 0\, .
	\end{align*}
	Similarly, dotting \eqref{gradientEq} by $\ve{n}$ and restricting the result to $\partial\Omega$ gives us the Neumann data
	\begin{align*}
		\left.\left[\left(\sum_{k=1}^{M}\frac{\alpha_k}{\rho_k}\right)\partial_{\vesmall{n}} p\right]\right|_{\partial\Omega} &= \left.\left[\sum_{k=1}^{M}\ve{n}\cdot\ve{f}_k  -\ve{n}\cdot\sum_{k=1}^{M}\partial_t(\alpha_k\ve{u}_k)\right]\right|_{\partial\Omega}\\
		&= \left.\left[\sum_{k=1}^{M}\ve{n}\cdot\ve{f}_k  -\sum_{k=1}^{M}\partial_t(\alpha_k\ve{n}\cdot\ve{u}_k)\right]\right|_{\partial\Omega} \\
		&= \sum_{k=1}^{M}\ve{n}\cdot\ve{f}_k\, ,
	\end{align*}
	since we are assuming $\ve{n}\cdot\ve{u}_k|_{\partial\Omega}=0$ for all $t$. Therefore, we find the initial pressure by solving the Neumann problem
	\begin{align}
		\nabla\cdot\left[\left(\sum_{k=1}^{M}\frac{\alpha_k^0}{\rho_k}\right)\nabla p^0\right] &= \nabla\cdot\sum_{k=1}^{M}\ve{f}_k^0 && \text{in} \ \, \Omega\, ,\\
		\partial_{\vesmall{n}} p^0 &= \left(\sum_{k=1}^{M}\frac{\alpha_k^0}{\rho_k}\right)^{-1}\sum_{k=1}^{M}\ve{n}\cdot\ve{f}_k^0 && \text{on} \ \, \partial\Omega\, ,
	\end{align}
	with $\ve{f}_k^0$ depending only on the initial data $(\alpha_k^0,\ve{u}_k^0)$. When using finite elements, this means simply finding $p^0\in Y_h\cap L^2_0(\Omega)$ such that
	\begin{equation*}
		\left\langle\nabla p^{0}, \left(\sum_{k=1}^{N}\frac{\alpha_k^{0}}{\rho_k}\right)\nabla q \right\rangle = \sum_{k=1}^{N}\left\langle\ve{f}_k^0, \nabla q \right\rangle
	\end{equation*}
	for all $q\in Y_h\cap L^2_0(\Omega)$.

	\section{Temporal stability analysis}\label{sec_analysis}
Some identities and inequalities will be useful for the analysis. Provided that $\ve{v}\cdot\ve{n}=0$ on $\partial\Omega$, the following equalities hold for any $(r,\ve{v},\ve{w})$ sufficiently regular \cite{Guermond2000}:
	\begin{align}
		\left\langle\ve{v}\cdot\nabla r + \frac{\nabla\cdot\ve{v}}{2}r,r\right\rangle &= 0\, ,
		\label{skewSymRho} \\
	\left\langle(r\nabla\ve{v})\ve{w} + \frac{\nabla\cdot(r\ve{w})}{2}\ve{v},\ve{v}\right\rangle &= 0\, .
		\label{skewSymU}
	\end{align}
Moreover, for any scalar-, vector- or tensor-valued pair ($\ve{a},\ve{b}$), there holds
	\begin{align}
		2\langle \ve{a}+\ve{b},\ve{a}
        \rangle &= \|\ve{a}\|^2 - \|\ve{b}\|^2 + \|\ve{a}+\ve{b}\|^2 \, . \label{identityBDF1}
	\end{align}
    
The analysis will also require the following Gronwall inequality \cite{Heywood1990}.
	\begin{lemma}[Discrete Gronwall inequality]\label{Lem:Gronwall-unconditional}
	Let $N\in\mathbb{N}$, and $\alpha,B,a_{n},b_{n},c_{n}$ be non-negative numbers for $n=1,\ldots,N$. Let us suppose that these numbers satisfy
\begin{align}\label{N-1}
    a_{N} + \sum_{n=1}^{N}b_n \leq B +  \alpha\sum_{n=1}^{N-1}a_n  \, .
\end{align}
Then, the following inequality holds:
\begin{align}
    a_{N} + \sum_{n=1}^{N}b_n \leq B\mathrm{e}^{\alpha N} \ \ \text{for} \ \, N\geq 1 \, .
	\end{align}
\end{lemma}

	We will next analyse the stability of our IMEX fractional-step scheme. Although we consider the discrete-in-time case without addressing spatial discretisation, the analysis also applies to conforming finite element spaces. As usual, we assume homogeneous Dirichlet conditions $\ve{u}_k^{n+1}|_{\partial\Omega}=\ve{0}$ and, for concision, consider $\ve{g}_k=\ve{0}$ for the analysis.
\begin{lemma}[Estimate on the explicit drag terms]\label{Lem:Drag} Assuming that all the volume fractions fulfil the lower bound
\begin{align}\label{Assumption-positivity}
    \alpha_k^n \geq \alpha_{\mathrm{min}} > 0  \ \text{a.e.~in}\ \Omega\, , \ \text{for all} \ \, n=0,\dots,N \ \text{and all} \ \,  k=1,\dots,M\, ,
\end{align}
and denoting as $\rho_{\mathrm{min}}$ the minimum density among $\rho_1,\dots,\rho_M$, the drag terms will satisfy
\begin{align}   &2\tau\sum_{k=1}^{M}\sum_{l=1}^{M}\left\langle\tilde{\gamma}_{kl}^{n}(\hat{\mathbf{u}}_k^{n} - \hat{\mathbf{u}}_l^{n}),\ve{u}_k^{n+1}\right\rangle \geq \tau\sum_{k=1}^{M}\sum_{l=1}^{M}\left\|\sqrt{\tilde{\gamma}_{kl}^{n}}(\hat{\mathbf{u}}_k^{n}-\hat{\mathbf{u}}_l^{n})\right\|^2 \nonumber  \\
    &-\sum_{k=1}^{M}\left\lbrace \rho_k\left\|\sqrt{\alpha_k^n} 
(\ve{u}_k^{n+1}-\hat{\mathbf{u}}_k^{n})\right\|^2 + [2\tau(M-1)\beta]^2\rho_k\left\|\sqrt{\alpha_k^n}\hat{\mathbf{u}}_k^n\right\|^2 \right\rbrace ,\label{dragEstimate}
\end{align}
where 
\begin{align}
      \beta := \max_{m=0,\dots,n}\left\lbrace\max_{k,l=1,\dots,M}\frac{1}{\sqrt{\rho_k\rho_l}}\left\|\frac{\tilde{\gamma}_{kl}^m}{\sqrt{\vphantom{\frac{\cdot}{\cdot}}\alpha_k^m\alpha_l^m}}\right\|_{\infty}\right\rbrace \leq \frac{D}{\rho_{\mathrm{min}}\alpha_{\mathrm{min}}} \, .\label{beta}  
\end{align}
\end{lemma}
\proof{Since $\tilde{\gamma}_{kl}^n = \tilde{\gamma}_{lk}^n$, we can write
\begin{align*}
&2\tau\sum_{k=1}^{M}\sum_{l=1}^{M}\left\langle\tilde{\gamma}_{kl}^{n}(\hat{\mathbf{u}}_k^{n} - \hat{\mathbf{u}}_l^{n}),\ve{u}_k^{n+1}\right\rangle\\
&= \tau\sum_{k=1}^{M}\sum_{l=1}^{M}\left\langle\tilde{\gamma}_{kl}^{n}(\hat{\mathbf{u}}_k^{n} - \hat{\mathbf{u}}_l^{n}),\ve{u}_k^{n+1}\right\rangle + \tau\sum_{l=1}^{M}\sum_{k=1}^{M}\left\langle\tilde{\gamma}_{lk}^{n}(\hat{\mathbf{u}}_l^{n} - \hat{\mathbf{u}}_k^{n}),\ve{u}_l^{n+1}\right\rangle\\
&= \tau\sum_{k=1}^{M}\sum_{l=1}^{M}\left\langle\tilde{\gamma}_{kl}^{n}(\hat{\mathbf{u}}_k^{n} - \hat{\mathbf{u}}_l^{n}),\ve{u}_k^{n+1}\right\rangle - \tau\sum_{k=1}^{M}\sum_{l=1}^{M}\left\langle\tilde{\gamma}_{kl}^{n}(\hat{\mathbf{u}}_k^{n} - \hat{\mathbf{u}}_l^{n}),\ve{u}_l^{n+1}\right\rangle\\
&= \tau\sum_{k=1}^{M}\sum_{l=1}^{M}\left\langle\tilde{\gamma}_{kl}^{n}(\hat{\mathbf{u}}_k^{n} - \hat{\mathbf{u}}_l^{n}),\ve{u}_k^{n+1}-\ve{u}_l^{n+1}\right\rangle \\
&= \tau\sum_{k=1}^{M}\sum_{l=1}^{M}\left\langle\tilde{\gamma}_{kl}^{n}(\hat{\mathbf{u}}_k^{n} - \hat{\mathbf{u}}_l^{n}),(\hat{\mathbf{u}}_k^{n}-\hat{\mathbf{u}}_l^{n}) + 
(\ve{u}_k^{n+1}-\hat{\mathbf{u}}_k^{n}) - (\ve{u}_l^{n+1}-\hat{\mathbf{u}}_l^{n})\right\rangle\\
&= \tau\sum_{k=1}^{M}\sum_{l=1}^{M}\left\|\sqrt{\tilde{\gamma}_{kl}^{n}}(\hat{\mathbf{u}}_k^{n}-\hat{\mathbf{u}}_l^{n})\right\|^2 + \left\langle\tilde{\gamma}_{kl}^{n}(\hat{\mathbf{u}}_k^{n} - \hat{\mathbf{u}}_l^{n}), 
(\ve{u}_k^{n+1}-\hat{\mathbf{u}}_k^{n}) - (\ve{u}_l^{n+1}-\hat{\mathbf{u}}_l^{n})\right\rangle\, .
\end{align*}
The first term on the right-hand side is non-negative, so it remains to estimate the other term. Recalling that $\tilde{\gamma}_{kk}=0$ and using Hölder's inequality, we can write
 \begin{align*}
    &\tau\left|\sum_{k=1}^{M}\sum_{l=1}^{M} \left\langle\tilde{\gamma}_{kl}^{n}(\hat{\mathbf{u}}_k^{n} - \hat{\mathbf{u}}_l^{n}), 
(\ve{u}_k^{n+1}-\hat{\mathbf{u}}_k^{n}) - (\ve{u}_l^{n+1}-\hat{\mathbf{u}}_l^{n})\right\rangle\right|\\
&= \tau\left|\sum_{k=1}^{M}\sum_{l=1}^{M}\left( \left\langle\tilde{\gamma}_{kl}^{n}(\hat{\mathbf{u}}_k^{n} - \hat{\mathbf{u}}_l^{n}), 
\ve{u}_k^{n+1}-\hat{\mathbf{u}}_k^{n}\right\rangle - \left\langle\tilde{\gamma}_{kl}^{n}(\hat{\mathbf{u}}_k^{n} - \hat{\mathbf{u}}_l^{n}), 
\ve{u}_l^{n+1}-\hat{\mathbf{u}}_l^{n}\right\rangle\right)\right|\\
&= 2\tau\left|\sum_{k=1}^{M}\sum_{l=1}^{M} \left\langle\tilde{\gamma}_{kl}^{n}(\hat{\mathbf{u}}_k^{n} - \hat{\mathbf{u}}_l^{n}), 
\ve{u}_k^{n+1}-\hat{\mathbf{u}}_k^{n}\right\rangle\right|\\
&= 2\tau\left|\sum_{k=1}^{M}\sum_{l=1}^{M} \left\langle\frac{\tilde{\gamma}_{kl}^{n}}{\rho_k\alpha_k^n}\sqrt{\rho_k\alpha_k^n}\,\hat{\mathbf{u}}_k^{n} - \frac{\tilde{\gamma}_{kl}^{n}}{\sqrt{\vphantom{\frac{\cdot}{\cdot}}\rho_k\alpha_k^n\rho_l\alpha_l^n}}\sqrt{\rho_l\alpha_l^n}\,\hat{\mathbf{u}}_l^{n},\sqrt{\rho_k\alpha_k^n} 
(\ve{u}_k^{n+1}-\hat{\mathbf{u}}_k^{n})\right\rangle\right|\\
&\leq 2\tau\sum_{k=1}^{M}\sum_{l=1}^{M}\left(\left\|\sqrt{\rho_k\alpha_k^n}\,\hat{\mathbf{u}}_k^n\right\|\left\|\frac{\tilde{\gamma}_{kl}^n}{\rho_k\alpha_k^n}\right\|_{\infty} + \left\|\sqrt{\rho_l\alpha_l^n}\,\hat{\mathbf{u}}_l^n\right\|\left\|\frac{\tilde{\gamma}_{kl}^n}{\sqrt{\vphantom{\frac{\cdot}{\cdot}}\rho_k\rho_l\alpha_k^n\alpha_l^n}}\right\|_{\infty}\right)\left\|\sqrt{\rho_k\alpha_k^n} 
(\ve{u}_k^{n+1}-\hat{\mathbf{u}}_k^{n})\right\|\\
&\leq 2\tau\beta\sum_{k=1}^{M}\sum_{\substack{l=1\\l\neq k}}^{M}\left(\left\|\sqrt{\rho_k\alpha_k^n}\,\hat{\mathbf{u}}_k^n\right\| + \left\|\sqrt{\rho_l\alpha_l^n}\,\hat{\mathbf{u}}_l^n\right\|\right)\left\|\sqrt{\rho_k\alpha_k^n} 
(\ve{u}_k^{n+1}-\hat{\mathbf{u}}_k^{n})\right\|,
 \end{align*}
introducing the positive constant $\beta$, defined in \eqref{beta}. Let us now use the short-hand notation $a_k=\left\|\sqrt{\vphantom{\frac{\cdot}{\cdot}}\rho_k\alpha_k^n}\,\hat{\mathbf{u}}_k^n\right\|$, $b_k = \left\|\sqrt{\vphantom{\frac{\cdot}{\cdot}}\rho_k\alpha_k^n} 
(\ve{u}_k^{n+1}-\hat{\mathbf{u}}_k^{n})\right\|$, so that Young's inequality gives us
\begin{align*}
    &2\tau\beta\sum_{k=1}^{M}\sum_{\substack{l=1\\l\neq k}}^{M}\left(\left\|\sqrt{\rho_k\alpha_k^n}\hat{\mathbf{u}}_k^n\right\| + \left\|\sqrt{\rho_l\alpha_l^n}\hat{\mathbf{u}}_l^n\right\|\right)\left\|\sqrt{\rho_k\alpha_k^n} 
(\ve{u}_k^{n+1}-\hat{\mathbf{u}}_k^{n})\right\| \\
&= 2\tau\beta\sum_{k=1}^{M}\sum_{\substack{l=1\\l\neq k}}^{M}(a_k b_k + a_l b_k)\\
&= 2\tau\beta\sum_{k=1}^{M}(M-1)a_k b_k + 2\sum_{k=1}^{M}\sum_{\substack{l=1\\l\neq k}}^{M} \tau\beta a_l b_k\\
&\leq \sum_{k=1}^{M}\left\lbrace\frac{(b_k)^2}{2} + 2[\tau\beta(M-1)a_k]^2\right\rbrace  + \sum_{k=1}^{M}\sum_{\substack{l=1\\l\neq k}}^{M}\left\lbrace\frac{(b_k)^2}{2(M-1)} + 2(M-1)(\tau\beta a_l)^2 \right\rbrace\\
&= \sum_{k=1}^{M}\left\lbrace\left[\frac{(b_k)^2}{2} + 2(\tau\beta(M-1)a_k)^2\right] + \left[(M-1)\frac{(b_k)^2}{2(M-1)} + (M-1)2(M-1)(\tau\beta a_k)^2\right]\right\rbrace \\
&= \sum_{k=1}^{M}\left\lbrace (b_k)^2 + (2\tau\beta(M-1)a_k)^2\right\rbrace \\
&= \sum_{k=1}^{M}\left\lbrace \rho_k\left\|\sqrt{\alpha_k^n} 
(\ve{u}_k^{n+1}-\hat{\mathbf{u}}_k^{n})\right\|^2 + [2\tau(M-1)\beta]^2\rho_k\left\|\sqrt{\alpha_k^n}\hat{\mathbf{u}}_k^n\right\|^2 \right\rbrace ,
\end{align*}
which when combined with the two previous estimates completes the proof.
}
	
	\begin{theorem}[Stability of velocities and pressure]
		Under the assumptions of Lemma \ref{Lem:Drag}, for any time-step size $\tau = T/N >0$ the fractional-step scheme \eqref{FS-momentum},\eqref{PPE},\eqref{projectedU} with $\ve{g}_k=\ve{0}$ ($k=1,\dots,M$) satisfies the stability estimate
		\begin{align}
&\sum_{k=1}^{M}\left(\rho_k\left\|\sqrt{\alpha_k^{N}}\hat{\mathbf{u}}_k^{N}\right\|^2 + \tau\mu_k\left\|\sqrt{\alpha_k^{N}}\,\nabla\ve{u}_k^{N}\right\|^2 + \frac{\tau^2}{\rho_k}\left\|\sqrt{\alpha_k^{N}}\,\nabla p^{N} \right\|^2 \right)   \nonumber\\
		&+\tau\sum_{n=1}^{N}\sum_{k=1}^{M}\left[\mu_k\left\|\sqrt{\vphantom{\alpha_k^{n+1}}\alpha_k^{n}}\,\nabla\ve{u}_k^{n}+\sqrt{\vphantom{\alpha_k^{n}}\alpha_k^{n-1}}\,\nabla^{\top}\ve{u}_k^{n-1}\right\|^2 + \sum_{l=1}^{M}\left\|\sqrt{\tilde{\gamma}_{kl}^{n-1}}(\hat{\mathbf{u}}_k^{n-1}-\hat{\mathbf{u}}_l^{n-1})\right\|^2\right]
			\nonumber\\
			&\leq \left(\sum_{k=1}^M B_k\right)\exp\left\lbrace{4\bigg[\frac{D(M-1)}{\rho_{\mathrm{min}}\alpha_{\mathrm{min}}}\bigg]^2 T\tau}\right\rbrace\, ,  \label{mainResult}
		\end{align}
where
\begin{align}
    B_k = \Bigg\lbrace{1+\bigg[\frac{2D(M-1)\tau}{\rho_{\mathrm{min}}\alpha_{\mathrm{min}}}\bigg]^2}\Bigg\rbrace\rho_k\left\|\sqrt{\alpha_k^{0}}\ve{u}_k^{0}\right\|^2 + \tau\mu_k\left\|\sqrt{\alpha_k^{0}}\,\nabla\ve{u}_k^{0}\right\|^2 + \frac{\tau^2}{\rho_k}\left\|\sqrt{\alpha_k^{0}}\,\nabla p^{0} \right\|^2.\label{Bk}  
\end{align}
\end{theorem}
	
\proof{We start by setting $\ve{v}=\ve{u}_k^{n+1}-\hat{\mathbf{u}}_k^{n+1}$ in Eq.~\eqref{weakUhat} and using the Cauchy--Schwarz inequality:
		\begin{align*}	\left\|\sqrt{\alpha_k^{n+1}}(\ve{u}_k^{n+1}-\hat{\mathbf{u}}_k^{n+1}) \right\|^2 &=  \frac{\tau}{\rho_k}\left\langle\sqrt{\vphantom{\alpha_k^{n+1}}\alpha_k^{n}\,}\,\nabla p^{n}-\sqrt{\alpha_k^{n+1}}\,\nabla p^{n+1},\sqrt{\alpha_k^{n+1}}(\ve{u}_k^{n+1}-\hat{\mathbf{u}}_k^{n+1})\right\rangle\\
			&\leq \frac{\tau}{\rho_k}\left\|\sqrt{\alpha_k^{n+1}}(\ve{u}_k^{n+1}-\hat{\mathbf{u}}_k^{n+1}) \right\| \left\|\sqrt{\vphantom{\alpha_k^{n+1}}\alpha_k^{n}\,}\,\nabla p^{n}-\sqrt{\alpha_k^{n+1}}\,\nabla p^{n+1}\right\|\, ,
		\end{align*}
		so that
		\begin{align}
			\sqrt{\rho_k}\left\|\sqrt{\alpha_k^{n+1}}(\ve{u}_k^{n+1}-\hat{\mathbf{u}}_k^{n+1}) \right\| \leq \frac{\tau}{\sqrt{\rho_k}}\left\|\sqrt{\alpha_k^{n+1}}\,\nabla p^{n+1} - \sqrt{\vphantom{\alpha_k^{n+1}}\alpha_k^{n}}\,\nabla p^{n}\right\| \, . \label{CS}
		\end{align}
		Taking $q=\tau p^{n+1}$ in \eqref{weakPPE} and using \eqref{divEndOfStep} gives
		\begin{align*}
			\sum_{k=1}^M\frac{\tau}{\rho_k}\left\langle\sqrt{\alpha_k^{n+1}}\,\nabla p^{n+1}-\sqrt{\vphantom{\alpha_k^{n+1}}\alpha_k^{n}}\,\nabla p^{n},\sqrt{\alpha_k^{n+1}}\nabla p^{n+1}\right\rangle &= \sum_{k=1}^M\left\langle-\nabla\cdot(\alpha_k^{n+1}\ve{u}_k^{n+1}),p^{n+1}\right\rangle\\
			&= \sum_{k=1}^M\left\langle\nabla\cdot[\alpha_k^{n+1}(\hat{\mathbf{u}}_k^{n+1}-\ve{u}_k^{n+1})],p^{n+1}\right\rangle\\
			&= \sum_{k=1}^M\left\langle\alpha_k^{n+1}(\ve{u}_k^{n+1}-\hat{\mathbf{u}}_k^{n+1}),\nabla p^{n+1}\right\rangle,
		\end{align*}
		where we have used integration by parts in the last step. Rewriting the left-hand side by using \eqref{identityBDF1} and \eqref{CS} yields
		\begin{align}
		&2\tau\sum_{k=1}^N\left\langle\alpha_k^{n+1}(\ve{u}_k^{n+1}-\hat{\mathbf{u}}_k^{n+1}),\nabla p^{n+1}\right\rangle \nonumber\\
			&=\sum_{k=1}^N\frac{ \tau^2}{\rho_k}\left(\left\|\sqrt{\alpha_k^{n+1}}\,\nabla p^{n+1} \right\|^2 - \left\|\sqrt{\vphantom{\alpha_k^{n+1}}\alpha_k^{n}}\,\nabla p_k^{n} \right\|^2 + \left\|\sqrt{\alpha_k^{n+1}}\,\nabla p^{n+1}-\sqrt{\vphantom{\alpha_k^{n+1}}\alpha_k^{n}}\,\nabla p^{n}\right\|^2\right)\nonumber\\
			&\geq \sum_{k=1}^N\left[\frac{\tau^2}{\rho_k}\left(\left\|\sqrt{\alpha_k^{n+1}}\,\nabla p^{n+1} \right\|^2 - \left\|\sqrt{\vphantom{\alpha_k^{n+1}}\alpha_k^{n}}\,\nabla p_k^{n} \right\|^2\right) + \rho_k\left\|\sqrt{\alpha_k^{n+1}}(\ve{u}_k^{n+1}-\hat{\mathbf{u}}_k^{n+1})\right\|^2\right].\label{estPPE}
		\end{align}
		Next, we set $\mathbf{v}=2\rho_k\hat{\mathbf{u}}_k^{n+1}$ in \eqref{weakUhat} to get
		\begin{align*}	&\rho_k\left(\left\|\sqrt{\alpha_k^{n+1}}\hat{\mathbf{u}}_k^{n+1} \right\|^2 -  \left\|\sqrt{\alpha_k^{n+1}}\ve{u}_k^{n+1} \right\|^2 +  \left\|\sqrt{\alpha_k^{n+1}}\left(\hat{\mathbf{u}}_k^{n+1}-\ve{u}_k^{n+1}\right) \right\|^2\right) \\
			&= 2\tau\left\langle\sqrt{\vphantom{\alpha_k^{n+1}}\alpha_k^{n}\,}\,\nabla p^{n},\sqrt{\alpha_k^{n+1}}\,\hat{\mathbf{u}}_k^{n+1}\right\rangle - 2\tau\left\langle\nabla p^{n+1},\alpha_k^{n+1}\hat{\mathbf{u}}_k^{n+1}\right\rangle .
		\end{align*}
		Adding over all $M$ phases gives us
		\begin{align}
			&\sum_{k=1}^{M}\rho_k\left(\left\|\sqrt{\alpha_k^{n+1}}\hat{\mathbf{u}}_k^{n+1} \right\|^2 -  \left\|\sqrt{\alpha_k^{n+1}}\ve{u}_k^{n+1} \right\|^2 +  \left\|\sqrt{\alpha_k^{n+1}}\left(\hat{\mathbf{u}}_k^{n+1}-\ve{u}_k^{n+1}\right) \right\|^2\right) \nonumber\\
			&= 2\tau\sum_{k=1}^{M}\left(\left\langle\sqrt{\vphantom{\alpha_k^{n+1}}\alpha_k^{n}\,}\,\nabla p^{n},\sqrt{\alpha_k^{n+1}}\,\hat{\mathbf{u}}_k^{n+1}\right\rangle - 2\tau\left\langle\nabla p^{n+1},\alpha_k^{n+1}\hat{\mathbf{u}}_k^{n+1}\right\rangle\right) \nonumber\\
			&= 2\tau\sum_{k=1}^{M}\left\langle\sqrt{\vphantom{\alpha_k^{n+1}}\alpha_k^{n}\,}\,\nabla p^{n},\sqrt{\alpha_k^{n+1}}\,\hat{\mathbf{u}}_k^{n+1}\right\rangle -2\tau\left\langle\nabla p^{n+1},\sum_{k=1}^{M}\alpha_k^{n+1}\hat{\mathbf{u}}_k^{n+1}\right\rangle \nonumber\\
			&= 2\tau\sum_{k=1}^{M}\left\langle\sqrt{\vphantom{\alpha_k^{n+1}}\alpha_k^{n}\,}\,\nabla p^{n},\sqrt{\alpha_k^{n+1}}\,\hat{\mathbf{u}}_k^{n+1}\right\rangle + 2\tau\left\langle p^{n+1},\nabla\cdot\sum_{k=1}^{M}\alpha_k^{n+1}\hat{\mathbf{u}}_k^{n+1}\right\rangle \nonumber\\		&=2\tau\sum_{k=1}^{M}\left\langle\sqrt{\vphantom{\alpha_k^{n+1}}\alpha_k^{n}\,}\,\nabla p^{n},\sqrt{\alpha_k^{n+1}}\,\hat{\mathbf{u}}_k^{n+1}\right\rangle , \label{estProjectedU}
		\end{align}
		again thanks to integration by parts and \eqref{divEndOfStep}. This, combined with estimate \eqref{estPPE}, yields 
		\begin{align} &\sum_{k=1}^{M}\rho_k\left(\left\|\sqrt{\alpha_k^{n+1}}\hat{\mathbf{u}}_k^{n+1} \right\|^2 -  \left\|\sqrt{\alpha_k^{n+1}}\ve{u}_k^{n+1} \right\|^2 +  2\left\|\sqrt{\alpha_k^{n+1}}\left(\hat{\mathbf{u}}_k^{n+1}-\ve{u}_k^{n+1}\right) \right\|^2\right) \nonumber\\
			&+ \sum_{k=1}^{M}\frac{\tau^2}{\rho_k}\left(\left\|\sqrt{\alpha_k^{n+1}}\nabla p^{n+1} \right\|^2 - \left\|\sqrt{\alpha_k^{n}\vphantom{\alpha_k^{n+1}}}\nabla p^{n} \right\|^2\right) \nonumber \\
			&\leq 2\tau\sum_{k=1}^{M}\left[ \left\langle\alpha_k^{n+1}(\ve{u}_k^{n+1}-\hat{\mathbf{u}}_k^{n+1}),\nabla p^{n+1}\right\rangle - \left\langle p^{n},\nabla\cdot\left(\sqrt{\alpha_k^{n+1}\alpha_k^{n}}\,\hat{\mathbf{u}}_k^{n+1}\right)\right\rangle\right].\label{estProjectionStep}
		\end{align}
		
		Let us now take $\ve{v} = 2\tau\ve{u}_k^{n+1}$ in \eqref{weakMomentum}, which gives us several terms:
		\begin{itemize}
			\item Acceleration terms:
			\begin{align}
				& \rho_k\left\langle (\alpha_k^{n+1}+\alpha_k^{n})\ve{u}_k^{n+1}-2\alpha_k^{n}\hat{\mathbf{u}}_k^{n},\ve{u}_k^{n+1}\right\rangle \nonumber\\
				&= \rho_k\left\|\sqrt{\alpha_k^{n+1}}\ve{u}_k^{n+1}\right\|^2 + \rho_k\left\langle \ve{u}_k^{n+1}-2\hat{\mathbf{u}}_k^{n},\alpha_k^{n}\ve{u}_k^{n+1}\right\rangle \nonumber \\ 
				&=\rho_k\left(\left\|\sqrt{\alpha_k^{n+1}}\ve{u}_k^{n+1}\right\|^2 -  \left\|\sqrt{\alpha_k^{n}\vphantom{\alpha_k^{n+1}}}\hat{\mathbf{u}}_k^{n}\right\|^2 + \left\|\sqrt{\alpha_k^{n}\vphantom{\alpha_k^{n+1}}}(\ve{u}_k^{n+1}-\hat{\mathbf{u}}_k^{n})\right\|^2\right)\, .\label{estAcc}
			\end{align}
			
			\item Convective terms:
			\begin{align}
				2\tau\left\langle \alpha_k^{n+1}\ve{u}_k^{n}\cdot\nabla\ve{u}_k^{n+1} + \frac{1}{2}[\nabla\cdot(\alpha_k^{n+1}\ve{u}_k^n)]\ve{u}_k^{n+1},\ve{u}_k^{n+1} \right\rangle = 0\, ,
			\end{align}
			due to the skew-symmetry result \eqref{skewSymU}.
			
			\item Viscous terms:
			\begin{align}
				&2\tau\mu_k\left\langle \sqrt{\alpha_k^{n+1}}\,\nabla\ve{u}_k^{n+1}+\sqrt{\vphantom{\alpha_k^{n+1}}\alpha_k^n}\,\nabla^{\top}\ve{u}_k^{n},\sqrt{\alpha_k^{n+1}}\,\nabla\ve{u}_k^{n+1}\right\rangle\nonumber \\
                &= \tau\mu_k\left(\left\|\sqrt{\alpha_k^{n+1}}\,\nabla\ve{u}_k^{n+1}\right\|^2 - \left\|\sqrt{\vphantom{\alpha_k^{n+1}}\alpha_k^n}\,\nabla^{\top}\ve{u}_k^{n}\right\|^2 + \left\|\sqrt{\alpha_k^{n+1}}\,\nabla\ve{u}_k^{n+1}+\sqrt{\vphantom{\alpha_k^{n+1}}\alpha_k^n}\,\nabla^{\top}\ve{u}_k^{n}\right\|^2 \right)\nonumber\\
                 &= \tau\mu_k\left(\left\|\sqrt{\alpha_k^{n+1}}\,\nabla\ve{u}_k^{n+1}\right\|^2 - \left\|\sqrt{\vphantom{\alpha_k^{n+1}}\alpha_k^n}\,\nabla\ve{u}_k^{n}\right\|^2 + \left\|\sqrt{\alpha_k^{n+1}}\,\nabla\ve{u}_k^{n+1}+\sqrt{\vphantom{\alpha_k^{n+1}}\alpha_k^n}\,\nabla^{\top}\ve{u}_k^{n}\right\|^2 \right).\label{viscousIMEX}
			\end{align}
			
			\item Drag terms: \begin{align}   &-2\tau\sum_{k=1}^{M}\sum_{l=1}^{M}\left\langle\tilde{\gamma}_{kl}^{n}(\hat{\mathbf{u}}_k^{n} - \hat{\mathbf{u}}_l^{n}),\ve{u}_k^{n+1}\right\rangle \leq -\tau\sum_{k=1}^{M}\sum_{l=1}^{M}\left\|\sqrt{\tilde{\gamma}_{kl}^{n}}(\hat{\mathbf{u}}_k^{n}-\hat{\mathbf{u}}_l^{n})\right\|^2 \nonumber  \\
    &+\sum_{k=1}^{M}\left\lbrace \rho_k\left\|\sqrt{\alpha_k^n} 
(\ve{u}_k^{n+1}-\hat{\mathbf{u}}_k^{n})\right\|^2 + [2\tau(M-1)\beta]^2\rho_k\left\|\sqrt{\alpha_k^n}\hat{\mathbf{u}}_k^n\right\|^2 \right\rbrace ,\label{dragEstimate2}
\end{align}
as proved in Lemma \ref{Lem:Drag}. 
\end{itemize}

Combining estimates \eqref{estAcc}--\eqref{dragEstimate2} yields, for the momentum equation,
		\begin{align*}
&\sum_{k=1}^{M}\left\lbrace\rho_k\left(\left\|\sqrt{\alpha_k^{n+1}}\ve{u}_k^{n+1}\right\|^2 -  \left\|\sqrt{\alpha_k^{n}\vphantom{\alpha_k^{n+1}}}\hat{\mathbf{u}}_k^{n}\right\|^2\right) + \tau\sum_{l=1}^{M}\left\|\sqrt{\tilde{\gamma}_{kl}^{n}}(\hat{\mathbf{u}}_k^{n}-\hat{\mathbf{u}}_l^{n})\right\|^2 +\right.   \\
&\left.\tau\mu_k\left(\left\|\sqrt{\alpha_k^{n+1}}\,\nabla\ve{u}_k^{n+1}\right\|^2 - \left\|\sqrt{\vphantom{\alpha_k^{n+1}}\alpha_k^n}\,\nabla^{\top}\ve{u}_k^{n}\right\|^2 + \left\|\sqrt{\alpha_k^{n+1}}\,\nabla\ve{u}_k^{n+1}+\sqrt{\vphantom{\alpha_k^{n+1}}\alpha_k^n}\,\nabla^{\top}\ve{u}_k^{n}\right\|^2 \right)\right\rbrace\\
&\leq \sum_{k=1}^{M}\left\lbrace-2\tau\left\langle\nabla p^{n},\sqrt{\alpha_k^{n+1}\alpha_k^{n}}\,\ve{u}_k^{n+1}\right\rangle + [2\tau(M-1)\beta]^2\rho_k\left\|\sqrt{\alpha_k^n}\hat{\mathbf{u}}_k^n\right\|^2\right\rbrace .
\end{align*}
Adding that to estimate \eqref{estProjectionStep} and defining $\alpha = [2\tau(M-1)\beta]^2$, we obtain
		\begin{align}	&\sum_{k=1}^{M}\rho_k\left(2\left\|\sqrt{\alpha_k^{n+1}}\left(\hat{\mathbf{u}}_k^{n+1}-\ve{u}_k^{n+1}\right) \right\|^2 + \left\|\sqrt{\alpha_k^{n+1}}\hat{\mathbf{u}}_k^{n+1}\right\|^2 -  \left\|\sqrt{\alpha_k^{n}\vphantom{\alpha_k^{n+1}}}\hat{\mathbf{u}}_k^{n}\right\|^2  \right) \nonumber\\
&+\tau\sum_{k=1}^{M}\mu_k\left(\left\|\sqrt{\alpha_k^{n+1}}\,\nabla\ve{u}_k^{n+1}\right\|^2 - \left\|\sqrt{\vphantom{\alpha_k^{n+1}}\alpha_k^n}\,\nabla\ve{u}_k^{n}\right\|^2 + \left\|\sqrt{\alpha_k^{n+1}}\,\nabla\ve{u}_k^{n+1}+\sqrt{\vphantom{\alpha_k^{n+1}}\alpha_k^n}\,\nabla^{\top}\ve{u}_k^{n}\right\|^2 \right)
			\nonumber\\
			&+ \sum_{k=1}^{M}\frac{\tau^2}{\rho_k}\left(\left\|\sqrt{\alpha_k^{n+1}}\nabla p^{n+1} \right\|^2 - \left\|\sqrt{\alpha_k^{n}\vphantom{\alpha_k^{n+1}}}\nabla p^{n} \right\|^2\right) + \tau\sum_{k=1}^{M}\sum_{l=1}^{M}\left\|\sqrt{\tilde{\gamma}_{kl}^{n}}(\hat{\mathbf{u}}_k^{n}-\hat{\mathbf{u}}_l^{n})\right\|^2 \nonumber \\
			&\leq \sum_{k=1}^{M}\left\lbrace 2\tau\left\langle\sqrt{\vphantom{\alpha_k^{n+1}}\alpha_k^{n}}\,\nabla p^{n}- \sqrt{\alpha_k^{n+1}}\,\nabla p^{n+1},\sqrt{\alpha_k^{n+1}}\,(\hat{\mathbf{u}}_k^{n+1}-\ve{u}_k^{n+1})\right\rangle + \alpha\rho_k\left\|\sqrt{\alpha_k^n}\hat{\mathbf{u}}_k^n\right\|^2 \right\rbrace \nonumber\\
			&= \sum_{k=1}^{M}\left\lbrace 2\rho_k\left\|\sqrt{\alpha_k^{n+1}}\left(\hat{\mathbf{u}}_k^{n+1}-\ve{u}_k^{n+1}\right) \right\|^2 + \alpha\rho_k\left\|\sqrt{\alpha_k^n}\hat{\mathbf{u}}_k^n\right\|^2 \right\rbrace ,\label{lastEstimate}
		\end{align}
		due to \eqref{weakUhat}. Then, the first term on the left-hand side cancels out the first one on the right-hand side. Rearranging some terms in \eqref{lastEstimate}, we have thus
\begin{align}	
&\sum_{k=1}^{M}\left[(\Psi_k^{n+1} - \Psi_k^{n}) + \tau\mu_k\left\|\sqrt{\alpha_k^{n+1}}\,\nabla\ve{u}_k^{n+1}+\sqrt{\vphantom{\alpha_k^{n+1}}\alpha_k^n}\,\nabla^{\top}\ve{u}_k^{n}\right\|^2 + \tau\sum_{l=1}^{M}\left\|\sqrt{\tilde{\gamma}_{kl}^{n}}(\hat{\mathbf{u}}_k^{n}-\hat{\mathbf{u}}_l^{n})\right\|^2 \right]\nonumber\\
			&\leq  \alpha\sum_{k=1}^{M}\rho_k\left\|\sqrt{\alpha_k^n}\hat{\mathbf{u}}_k^n\right\|^2 ,\label{organised}
		\end{align}
in which $ \Psi_k^{n+1} :=  \rho_k\Big\|\sqrt{\alpha_k^{n+1}}\hat{\mathbf{u}}_k^{n+1}\Big\|^2 + \tau\mu_k\Big\|\sqrt{\alpha_k^{n+1}}\,\nabla\ve{u}_k^{n+1}\Big\|^2 + \frac{\tau^2}{\rho_k}\Big\|\sqrt{\alpha_k^{n+1}}\nabla p^{n+1} \Big\|^2 \, .$

\noindent Adding from $n=0$ to $n=N-1$ gives
\begin{align*}	
&\sum_{k=1}^{M}\left\lbrace\Psi_k^{N}  + \tau\sum_{n=1}^{N}\left[\mu_k\left\|\sqrt{\vphantom{\alpha_k^{n+1}}\alpha_k^{n}}\,\nabla\ve{u}_k^{n}+\sqrt{\vphantom{\alpha_k^{n}}\alpha_k^{n-1}}\,\nabla^{\top}\ve{u}_k^{n-1}\right\|^2 + \sum_{l=1}^{M}\left\|\sqrt{\tilde{\gamma}_{kl}^{n}}(\hat{\mathbf{u}}_k^{n}-\hat{\mathbf{u}}_l^{n})\right\|^2\right]  \right\rbrace\\
&  \leq  \sum_{k=1}^{M}\left(\Psi_k^0 + \alpha\sum_{n=0}^{N-1}\rho_k\left\|\sqrt{\alpha_k^n}\hat{\mathbf{u}}_k^n\right\|^2\right)\\
&  =  \sum_{k=1}^{M}\left(B_k + \alpha\sum_{n=1}^{N-1}\rho_k\left\|\sqrt{\alpha_k^n}\hat{\mathbf{u}}_k^n\right\|^2\right)\\
&  \leq  \sum_{k=1}^{M}\left(B_k  + \alpha\sum_{n=1}^{N-1}\Psi_k^n\right),
\end{align*}
with $B_k$ defined in \eqref{Bk}. The proof is concluded by using the Gronwall Lemma \ref{Lem:Gronwall-unconditional} with 
\begin{align*}
a_n &= \sum_{k=1}^{M}\Psi_k^n \, , \\
b_n &= \tau\sum_{k=1}^{M}\bigg(\mu_k\bigg\|\sqrt{\vphantom{\alpha_k^{n+1}}\alpha_k^{n}}\,\nabla\ve{u}_k^{n}+\sqrt{\vphantom{\alpha_k^{n}}\alpha_k^{n-1}}\,\nabla^{\top}\ve{u}_k^{n-1}\bigg\|^2 + \sum_{l=1}^{M}\bigg\|\sqrt{\tilde{\gamma}_{kl}^{n-1}}(\hat{\mathbf{u}}_k^{n-1}-\hat{\mathbf{u}}_l^{n-1})\bigg\|^2\bigg).
\end{align*}
}
\begin{remark}
    For $\tau \lesssim \frac{\rho_{\mathrm{min}}\alpha_{\mathrm{min}}}{(M-1)D}$, it is possible to prove stability without the exponential factor seen in \eqref{mainResult}. That factor is in any case a mild one, as it tends to $1$ as $\tau\rightarrow 0$.
\end{remark}	

A crucial feature of the stability result \eqref{mainResult} is that, as just seen, it is proved without invoking the volume-fraction equations or their stability, which is only possible due to the reformulation we used for the convective term. This leaves us free to discretise (and even stabilise) the volume fractions in a variety of ways.

\section{Different formulations for the convection problem}\label{sec_convection}
Let us drop here the phase index $k$, as the following discussion applies very generally to the convection equation
\begin{equation}
    \partial_t\alpha + \nabla\cdot(\alpha\ve{u}) = 0\, ,\label{convectionGeneric}
\end{equation}
where we consider again non-penetrable boundaries ($\ve{n}\cdot\ve{u}|_{\partial\Omega} = 0$). Various challenges arise when $\ve{u}$ is not solenoidal; an important one is that, even if the initial condition $\alpha(0)$ is bounded in some interval $[a,b]$, there is no guarantee that $\alpha(t)$ will obey this constraint---even at the continuous level \cite{Kuzmin2012}. This shows that bound preservation is generally also a \textsl{modelling} matter, even before becoming a numerical one. Indeed, multiplying \eqref{convectionGeneric} by $2\alpha$ and integrating over $\Omega$ gives
\begin{align*}
    0 &= \int_{\Omega}2\alpha\partial_t\alpha \, \mathrm{d}\Omega + 2\int_{\Omega}\alpha\nabla\cdot(\alpha\ve{u}) \, \mathrm{d}\Omega\\
    &= \int_{\Omega}\partial_t(\alpha^2) \, \mathrm{d}\Omega + \int_{\Omega}\alpha\left[2\ve{u}\cdot\nabla\alpha +
    (\nabla\cdot\ve{u})\alpha + (\nabla\cdot\ve{u})\alpha\right]\mathrm{d}\Omega\\
    &= \frac{\mathrm{d}}{\mathrm{d}t}\|\alpha\|^2 + \int_{\Omega}(\nabla\cdot\ve{u})\alpha^2\,\mathrm{d}\Omega\, ,
\end{align*}
where we have used \eqref{skewSymRho}. Therefore, by defining the positive and negative parts of a scalar field $\phi$ as $\phi^+ := \max\lbrace\phi(\ve{x}),0 \rbrace$ and $\phi^-:=\phi - \phi^+$, it follows that
\begin{align}
   \frac{\mathrm{d}}{\mathrm{d}t}\|\alpha\|^2 + \left\|\sqrt{(\nabla\cdot\ve{u})^+}\,\alpha\right\|^2  =  \left\|\sqrt{-(\nabla\cdot\ve{u})^-}\,\alpha\right\|^2 \, .
   \label{instabilityAlpha}
\end{align}
It is easy to construct examples where $\nabla\cdot\ve{u} < 0$, in which case \eqref{instabilityAlpha} implies that $\|\alpha\|$ will increase indefinitely with time.  This property, of course, translates to the discrete level as numerical instability. Many approaches exist to overcome this, such as overshoot limiters and flux correction---we refer the reader to two recent books on this topic \cite{Kuzmin2023,Barrenechea2025}.

Motivated by this discussion, we can introduce a changed variable $\varphi := \sqrt{\alpha}$, so that substituting $\alpha=\varphi^2$ into Eq.~\eqref{mass} yields
		\begin{align*}
			0 &=  \partial_t\alpha + \ve{u}\cdot\nabla\alpha + (\nabla\cdot\ve{u})\alpha\\
			&= \partial_t(\varphi^2) + \ve{u}\cdot\nabla(\varphi^2) + (\nabla\cdot\ve{u})\varphi^2\\
			&= 2\varphi\left(\partial_t\varphi + \ve{u}\cdot\nabla\varphi +  \frac{\nabla\cdot\ve{u}}{2}\varphi\right),
		\end{align*}
which gives us the new equation
\begin{equation}
    \partial_t\varphi + \ve{u}\cdot\nabla\varphi +  \frac{\nabla\cdot\ve{u}}{2}\varphi = 0\, .\label{mass2}
\end{equation}
This leads to
\begin{align*}
0&=\int_{\Omega}2\varphi\Big(\partial_t\varphi + \ve{u}\cdot\nabla\varphi +  \frac{\nabla\cdot\ve{u}}{2}\varphi\Big)\mathrm{d}\Omega \\
&=\int_{\Omega}\partial_t(\varphi^2)\,\mathrm{d}\Omega + \int_{\Omega}\left[\ve{u}\cdot\nabla(\varphi^2) +  \varphi^2\nabla\cdot\ve{u}\right]\mathrm{d}\Omega \\
&= \frac{\mathrm{d}}{\mathrm{d}t}\int_{\Omega}\varphi^2\,\mathrm{d}\Omega + \int_{\Omega} \nabla\cdot(\varphi^2\ve{u})\,\mathrm{d}\Omega\\
&=\frac{\mathrm{d}}{\mathrm{d}t}\|\varphi\|^2 + \int_{\partial\Omega} \varphi^2\ve{n}\cdot\ve{u}\,\mathrm{d}\Gamma\\
&= \frac{\mathrm{d}}{\mathrm{d}t}\|\varphi\|^2\, ,
\end{align*}
that is, the $L^2(\Omega)$ norm of $\varphi$ is conserved over time. It can be shown that this stability is inherited at the discrete level by various time-stepping schemes, such as Crank--Nicolson, BDF2, implicit Euler, and the IMEX variant considered here:
\begin{align}
    \frac{\varphi^{n+1}-\varphi^{n}}{\tau} + \ve{u}^n\cdot\nabla\varphi^{n+1} + \frac{\nabla\cdot\ve{u}^n}{2}\varphi^{n+1} = 0\, .\label{varphiBDF1}
\end{align}

Here, and generally in numerical methods for conservation laws, there are trade-offs. We can choose to either tackle \eqref{convectionGeneric} directly---and be mass-conserving but potentially unstable---or discretise the reformulated equation \eqref{mass2} instead---and be unconditionally stable but only approximately conservative (mass conversation will follow the numerical convergence of the scheme). 

Besides mass conservation, another important constraint in our particular application is that $\alpha$ must be bounded between $0$ and $1$. Various changes of variables can be introduced to (en)force such bounds \cite{Hassler2020,Dzanic2023}, and the simple one we propose here is
\begin{equation}
    \alpha(\phi) = \left(\frac{\phi}{1+|\phi|}\right)^2 \in [0,1) \quad \text{for all} \ \, \phi\in\mathbb{R}\, ,
    \label{newVar}
\end{equation}
where
\begin{equation}
    \phi = \frac{\sqrt{\alpha}}{1-\sqrt{\alpha}}\, .
\end{equation}
The reader can verify that plugging \eqref{newVar} into \eqref{convectionGeneric} and using the chain rule produces the nonlinear equation
\begin{equation}
    \partial_t\phi + \ve{u}\cdot\nabla\phi +  \frac{\nabla\cdot\ve{u}}{2}\phi + \frac{\nabla\cdot\ve{u}}{2}|\phi|\phi = 0\, .\label{mass3}
\end{equation}
To circumvent the nonlinearity, we use the IMEX discretisation
\begin{equation}
    \frac{\phi^{n+1}-\phi^{n}}{\tau} + \ve{u}^{n}\cdot\nabla\phi^{n+1} +  \frac{\nabla\cdot\ve{u}^n}{2}\phi^{n+1} + \frac{\nabla\cdot\ve{u}^n}{2}|\phi^n|\phi^{n+1} = 0\, .\label{newVarIMEX}
\end{equation}
Although we shall not analyse this formulation, we expect once again a trade-off: the resulting volume fraction $\alpha^{n+1} = \alpha(\phi^{n+1})$ will be pointwise bounded in $[0,1)$, but probably at the cost of some time-step restriction for stability. In fact, most bound-preserving techniques induce CFL-like conditions \cite{Kuzmin2012,Dzanic2023}.

While both changes of variable presented here will be numerically tested, we shall only analyse the unconditionally stable one \eqref{varphiBDF1} among them. Now back to the multiphase context, we define another 
$H^1$-conforming finite element space $Z_h$. Then, prior to the velocity-pressure steps we update each volume fraction by first solving the weak problems to find $\varphi_k^{n+1}\in Z_h$ such that 
		\begin{equation}
			\left\langle\frac{\varphi_k^{n+1} - \varphi_k^n}{\tau} + \ve{u}_k^n\cdot\nabla\varphi_k^{n+1} + \frac{\nabla\cdot\ve{u}_k^{n}}{2}\varphi_k^{n+1},\frac{\zeta}{\tau} + \Big(\ve{u}_k^n\cdot\nabla\zeta + \frac{\nabla\cdot\ve{u}_k^{n}}{2}\zeta\Big)\chi\right\rangle = 0  \label{LS}
		\end{equation}
		for all $\zeta\in Z_h$, where $\chi=0$ indicates the standard Galerkin method, whereas $\chi=1$ leads to the least-squares Galerkin formulation, which is symmetric and usually more stable in the presence of sharp gradients \cite{Ern2004}. After updating $\varphi_k$, we simply evaluate $\alpha_k^{n+1} = (\varphi_k^{n+1})^2$. We will now prove the following stability result.
        
        \begin{theorem}[Stability of the volume fractions] Assuming $\ve{n}\cdot\ve{u}_k|_{\partial\Omega}=0$ at all times, scheme \eqref{LS} yields, for any time-step size $\tau > 0$, 
		\begin{align}
			\|\varphi_k^{N}\|^2 + \sum_{n=1}^{N}\left(\|\delta\varphi_k^{n} + \chi\tau\psi_k^{n} \|^2 + \chi\tau^2\|\psi_k^{n}\|^2\right) = \|\varphi_k^{0}\|^2\, ,
			\label{stabilityPhiLS}
		\end{align} 
		where $\delta\varphi_k^{n+1}:=\varphi_k^{n+1}-\varphi_k^{n}$ and $\psi_k^{n+1}:=\ve{u}_k^n\cdot\nabla\varphi_k^{n+1}+\frac{1}{2}\nabla\cdot\ve{u}_k^{n+1}\varphi_k^{n+1}$. This implies
		\begin{align}
			\|\alpha_k^{N}\|_{L^1(\Omega)} \leq \|\alpha_k^{0}\|_{L^1(\Omega)} \, ,
		\end{align}
		since $\alpha_k^{N} = (\varphi_k^{N})^2$.
	\end{theorem}
	\proof{We set $\zeta=2\tau^2\varphi_k^{n+1}$ in \eqref{LS} to get
		\begin{align*}
			0 &= 2\left\langle \delta\varphi_k^{n+1} + \tau\psi_k^{n+1},\varphi_k^{n+1} + \chi\tau\psi_k^{n+1}\right\rangle\\
			&=2\left\langle \delta\varphi_k^{n+1},\varphi_k^{n+1}\right\rangle
			+2\chi\left\langle \delta\varphi_k^{n+1}, \tau\psi_k^{n+1}\right\rangle + 2\chi\|\tau\psi_k^{n+1}\|^2
			+ 2\tau\left\langle  \psi_k^{n+1},\varphi_k^{n+1}\right\rangle\\
            &=\|\varphi_k^{n+1}\|^2 - \|\varphi_k^{n}\|^2+\|\delta\varphi_k^{n+1}\|^2
			+2\chi\langle \delta\varphi_k^{n+1}, \tau\psi_k^{n+1}\rangle + 2\chi\|\tau\psi_k^{n+1}\|^2
			+ 2\tau\langle  \psi_k^{n+1},\varphi_k^{n+1}\rangle\, ,
            \end{align*}
          where we have used \eqref{identityBDF1}. The last term on the right-hand side is zero due to identity \eqref{skewSymRho}. Moreover, since $\chi\in\lbrace 0,1\rbrace$, we can write $\chi = \chi^2$. Hence:  
            \begin{align*}
			0 &= \|\varphi_k^{n+1}\|^2 - \|\varphi_k^{n}\|^2+\|\delta\varphi_k^{n+1}\|^2
			+2\chi\left\langle \delta\varphi_k^{n+1}, \tau\psi_k^{n+1}\right\rangle + 2\chi^2\|\tau\psi_k^{n+1}\|^2\\
			&=  \|\varphi_k^{n+1}\|^2 - \|\varphi_k^{n}\|^2 + \|\delta\varphi_k^{n+1} + \chi\tau\psi_k^{n+1} \|^2 + \chi^2\tau^2\|\psi_k^{n+1}\|^2\, ,
		\end{align*}
which when added up from $n=0$ to $n=N-1$ completes the proof.}
	
\section{Numerical examples}\label{sec_examples}
	In this section, we assess the accuracy and the stability of our IMEX fractional-step scheme. All the tests were computed using quadrilateral Lagrangian finite elements. Unless where otherwise stated, (1) second- and first-order interpolation are used for velocities and pressure, respectively, and (2) the volume fractions are evolved using the unconditionally stable formulation \eqref{LS}, with the changed variables $\varphi_k^{n+1}=\sqrt{\alpha_k^{n+1}}$ discretised using first-order elements. Only two-phase ($M=2$) examples are solved, and the velocity boundary conditions are such that $\ve{n}\cdot\ve{u}_k=0$ on $\partial\Omega$ at all times, which fits our theory. Since there is no inflow through $\partial\Omega$, there are no boundary conditions for the volume fractions. The equations for the changed variables $\varphi_k$ or $\phi_k$ are formulated in the least-squares variant ($\chi=1$ in \eqref{LS}). Since the pressure is only defined (at each time) up to a constant, we employ the usual zero-mean scaling 
    \begin{align*}
        \int_{\Omega} p \, \mathrm{d}\Omega  = 0\quad\text{for all} \ \, t\, .
    \end{align*}

    Whenever the exact solutions $(p,\ve{u}_1,\ve{u}_2)$ are known, we will compute the relative errors at the last time step $N$ as 
    \begin{align}
        &e_p := \frac{\|p^N - p(T) \|}{\| p(T)\|}\, ,\label{errorP}\\
        &e_{\vesmall{u}} := \frac{\left\|\nabla\ve{u}_r^N - \nabla\ve{u}_r(T) \right\|}{\| \nabla\ve{u}_r(T)\|}\, , \quad \text{where} \ \, \ve{u}_r := \ve{u}_2-\ve{u}_1\, .
    \end{align}
    We can also measure the errors
    \begin{align}
        &e_{\text{div}} := \frac{1}{|\Omega|}\left\|\nabla\cdot(\alpha_1\hat{\mathbf{u}}_1^N + \alpha_2\hat{\mathbf{u}}_2^N) \right\|\, ,\\
        &e_{\alpha} := \frac{\left|\|\varphi_1^N\|^2 +\|\varphi_2^N\|^2 - \|1\|^2\right| }{\| 1\|^2} = \frac{1}{|\Omega|}\|\alpha_1^N +\alpha_2^N - 1 \|_{L^1(\Omega)}\, ,\label{errorDiv}
    \end{align}
    where $|\Omega|$ denotes (in 2D) the area of $\Omega$.
	
	\subsection{Convergence test with linear drag}
\begin{figure}[ht!]
		\centering
		\includegraphics[trim = 0 10 0 20,clip, width = .75\textwidth]{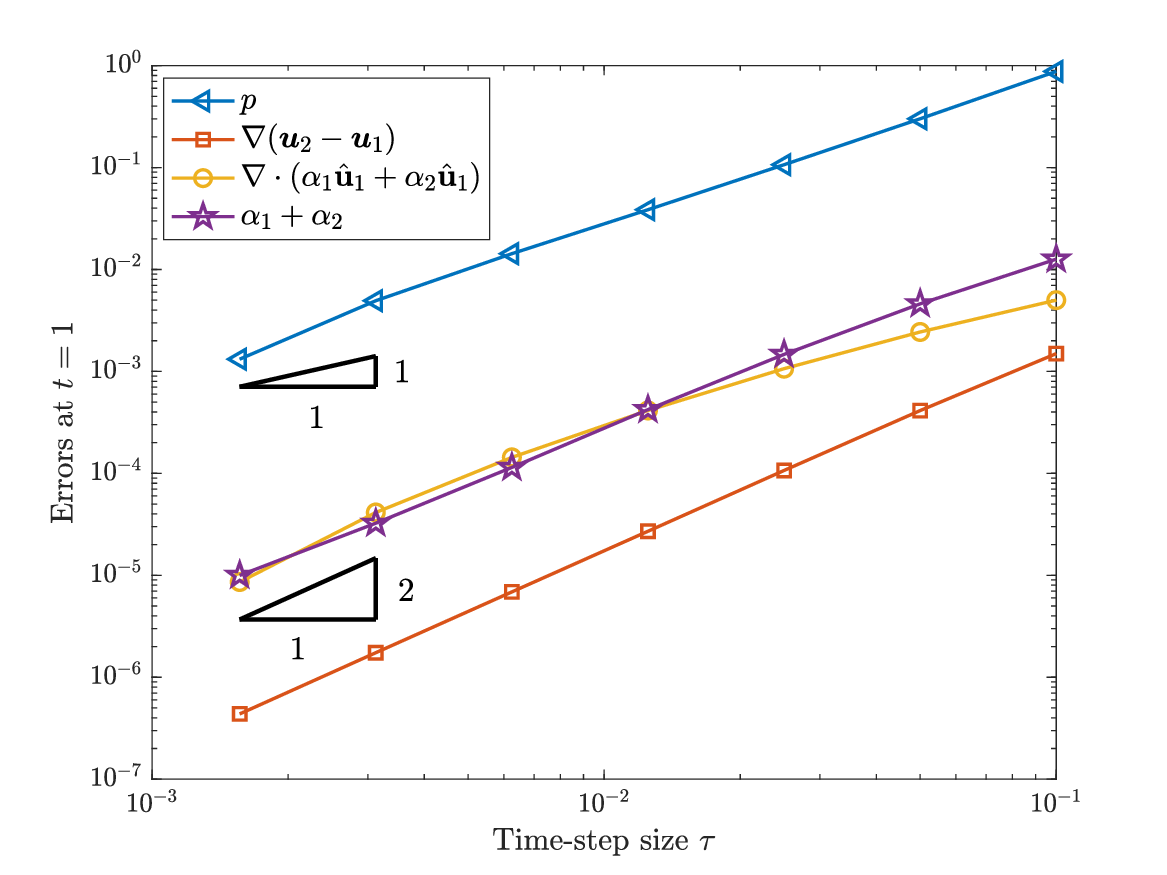}
		\caption{Temporal convergence test for a problem with constant volume fractions and linear drag. The relative errors \eqref{errorP}--\eqref{errorDiv} reveal a numerical order of convergence of at least $1$ for all four quantities.}
		\label{eoc1}
	\end{figure}

	The first numerical test considers a linear drag model with $\gamma_{12}=[4(t+1)]^{-1}$. The domain is the unit circle centred at $(0,0)$, there are no forcing terms ($\ve{g}_1=\ve{g}_2=\ve{0}$), and the fluid parameters are $\mu_1=\mu_2=\rho_1=\rho_2=1$. With that, we consider the solution
	\begin{align*}
		\ve{u}_1  &= f(t)
		\begin{pmatrix}
			-y\\
			x
		\end{pmatrix}, \ \
		\ve{u}_2  = -\ve{u}_1\,, \ \ p=f^2(t)\left(\frac{x^2+y^2}{2}-\frac{1}{4}\right), \ \ \alpha_1 = \alpha_2 \equiv \frac{1}{2}\, ,
	\end{align*}
	where $f(t)=(1+t)^{-1}$; the initial and boundary data are computed from the analytical expressions. The mesh has 12,288 elements, and the convergence study starts with $\tau=0.1$, with seven temporal refinements ($\tau\rightarrow\tau/2$) then applied. The errors defined in \eqref{errorP}--\eqref{errorDiv}, for $T=1$, are shown in Figure \ref{eoc1}. All quantities converge at least linearly, and some superconvergence is observed, which is common for incremental pressure-correction methods \cite{Barrenechea2024}.

\subsection{Convergence test with nonlinear drag}
\begin{figure}[ht!]
		\centering
		\includegraphics[trim = 0 10 0 20,clip, width = .75\textwidth]{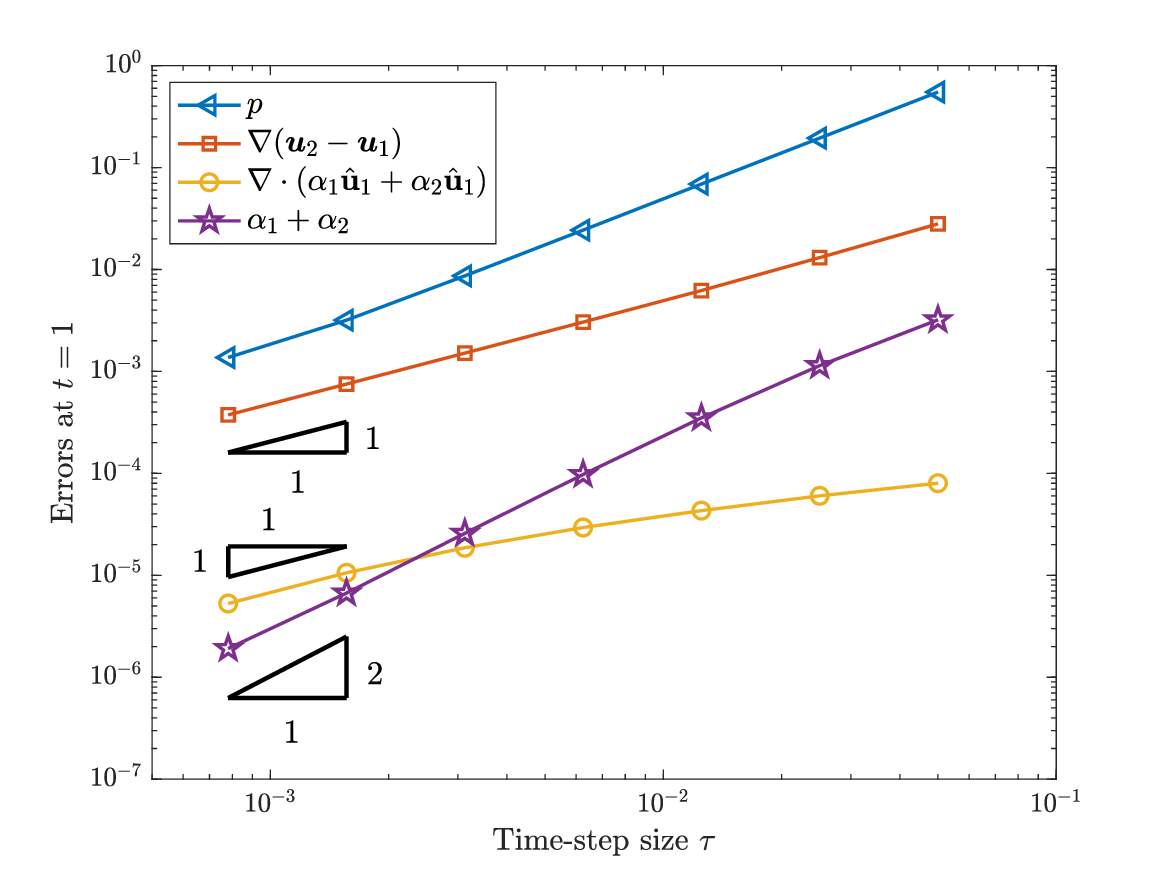}
		\caption{Temporal convergence test for a problem with quadratic drag. The relative errors at $t=1$ confirm at least linear convergence for all four quantities.}
		\label{eoc2}
	\end{figure}

The next test considers quadratic drag, $\gamma_{12} = 4|\ve{u}_1-\ve{u}_2|$, and the forcing terms
\begin{align*}
        \ve{g}_1 = \ve{0}\, , \ \ \ve{g}_2 = \frac{1}{4(1+t)^2}\frac{2-\sqrt{x^2+y^2}}{1-\sqrt{x^2+y^2}}\begin{pmatrix}
				    y \\ -x
				\end{pmatrix}.
    \end{align*}
The fluid parameters are $\mu_1=\rho_1=1$ and $\mu_2=\rho_2=4$. With that and the corresponding boundary and initial data, the solution is
	\begin{align*}
		\ve{u}_1  &= f(t)
		\begin{pmatrix}
			-y\\
			x
		\end{pmatrix}, \ \
		\ve{u}_2  = \frac{1}{2}\ve{u}_1\,, \ \ p=f^2(t)\left(\frac{x^2+y^2}{2}-\frac{5}{32}\right), \ \ \alpha_1 = 1-\alpha_2 = \sqrt{x^2+y^2}\, ,
	\end{align*}
	again with $f(t) = (1+t)^{-1}$. The domain is an annulus with inner and outer radii equal to $\nicefrac{1}{4}$ and $\nicefrac{3}{4}$, respectively, centred at the origin. The mesh has 24,576 elements with linear and quadratic interpolation for velocity and pressure, respectively. The temporal refinement is similar as before, this time starting with $\tau=0.05$. The results shown in Figure \ref{eoc2} confirm again at least linear convergence for all four quantities.

\subsection{Dispersed Rayleigh--Taylor flow}
\begin{figure}[ht!]
		\centering
		\includegraphics[trim = 80 90 10 50,clip, width = .8\textwidth]{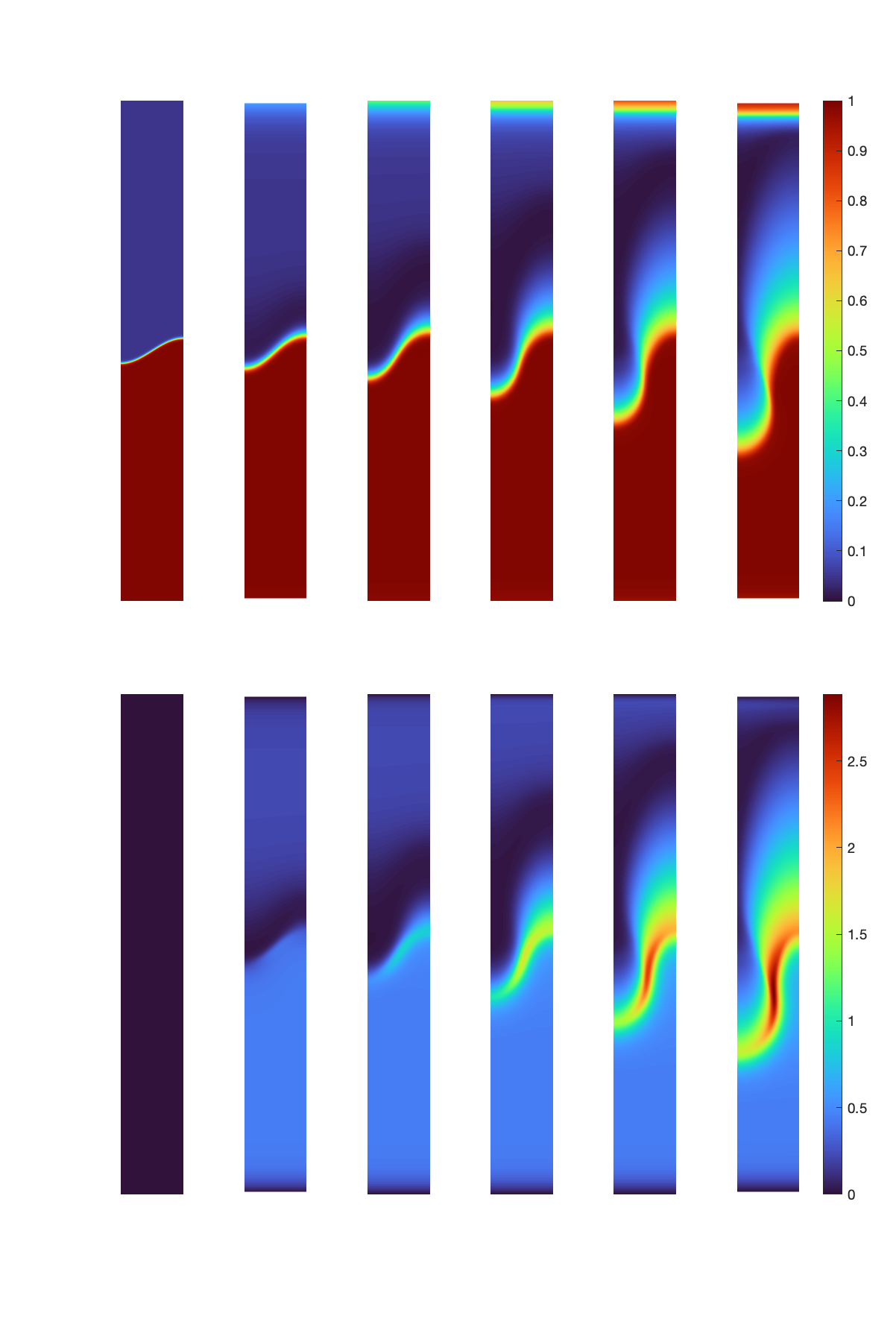}
		\caption{Dispersed Rayleigh--Taylor flow: distribution of $\alpha_2$ (top) and $\gamma_{12}$ (bottom) at $t=0,1,2,3,4,5$ (from left to right). The color scale for $\gamma_{12}$ goes from 0 to $\gamma_{12}^{\mathrm{max}} \approx 2.8838$.}
		\label{alphaAndDrag}
	\end{figure}

	Our final experiment is inspired by the Rayleigh--Taylor instability, which is a popular benchmark for (separated) two-phase flows \cite{Guermond2000}. The standard setup considers two fluids initially at rest in $\Omega = (0,0.5)\times(-2,2)$, with the denser one on top. To adapt this problem to the dispersed setting, we consider the initial phase distribution as
	\begin{align*}
		\alpha_2^0 = \frac{0.99+0.05}{2} + \frac{0.99-0.05}{2}\tanh\left(40y + 4\cos 2\pi x\right), \ \ \alpha_1^0 = 1-\alpha_2^0\, ,
	\end{align*}
	which means the lighter fluid ($k=2$) occupies $99\%$ of the lower part of $\Omega$ and only $5\%$ of the upper part, see Figure \ref{alphaAndDrag} (left). The lateral walls are free-slip boundaries, while the top and bottom ones have no slip. The system is under gravity $\ve{g}_1=\ve{g}_2=(0,-1)^{\top}$, and the fluid parameters are $\rho_1=1$, $\rho_2=3$, $\mu_1=0.1$, $\mu_2=0.3$. An Euler-type drag is considered: $\gamma_{12} = 10\alpha_2|\ve{u}_2-\ve{u}_1|$. The spatial mesh contains 80,000 square elements, and the time-step size is $\tau=0.005$. The volume fractions are evolved through the bound-preserving formulation \eqref{newVarIMEX}, with the changed variables $\phi_k$ discretised quadratically in space.

The distribution of $\alpha_2$ and $\gamma_{12}$ at different times is shown in Figure \ref{alphaAndDrag}. As in the classical separated two-phase benchmark, a mushroom-like pattern forms over time. However, the dispersed nature of the present problem allows some of the lighter fluid to rise to the top and start collecting underneath the upper wall, which cannot happen in the separated case. 
	
\begin{figure}[ht!]
		\centering
		\includegraphics[trim = 50 50 95 40,clip, width = 1\textwidth]{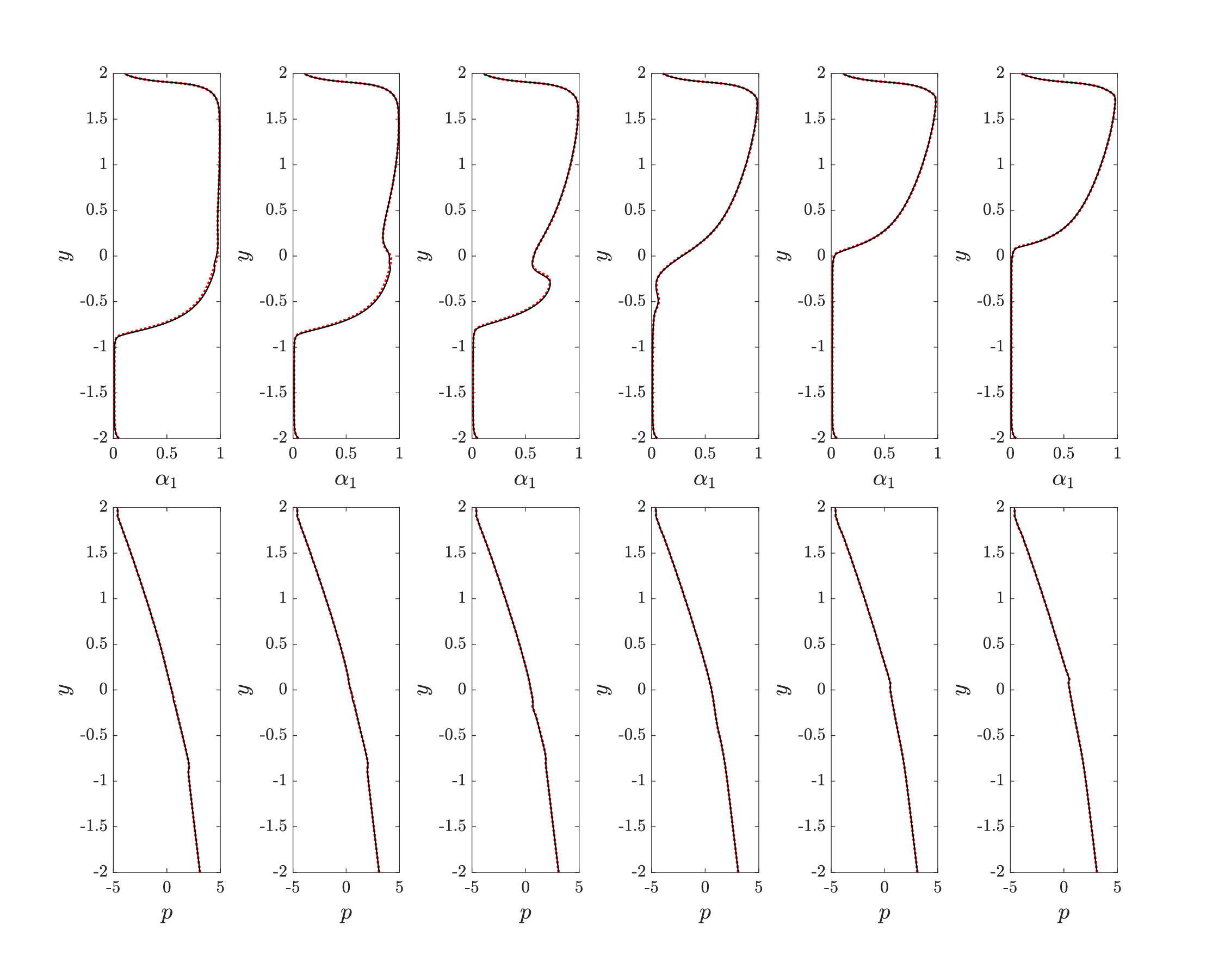}
		\caption{Dispersed Rayleigh--Taylor flow: $\alpha_1$ (top) and $p$ (bottom) profiles at $t=5$ for $x=0,0.1,0.2,0.3,0.4,0.5$ (from left to right). The fractional-step (full lines) and monolithic (dashed lines) solutions are in very good agreement.}
		\label{alpha1}
	\end{figure}

Since multiphase benchmark solutions are difficult to find, we compare our numerical results to the solution obtained through a (semi-)monolithic scheme with fully implicit viscous terms and IMEX drag forces:
\begin{align*}
\frac{\varphi_1^{n+1} - \varphi_1^n}{\tau} + \ve{u}_1^n\cdot\nabla\varphi_1^{n+1} + \frac{\nabla\cdot\ve{u}_1^{n}}{2}\varphi_1^{n+1} &= 0\, ,\\
\frac{\varphi_2^{n+1} - \varphi_2^n}{\tau} + \ve{u}_2^n\cdot\nabla\varphi_2^{n+1} + \frac{\nabla\cdot\ve{u}_2^{n}}{2}\varphi_2^{n+1} &= 0\, ,\\
\rho_1\bigg[\frac{\alpha_1^{n+1}+\alpha_1^{n}}{2\tau}\ve{u}_1^{n+1}   + \alpha_1^{n+1}\ve{u}_1^n\cdot\nabla\ve{u}_1^{n+1} + \frac{1}{2}\nabla\cdot(\alpha_1^{n+1}\ve{u}_1^n)\ve{u}_1^{n+1}\bigg]  &-\nabla\cdot(2\mu_1\alpha^{n+1}_1\nabla^{\mathrm{s}}\ve{u}^{n+1}_1)  \\
		+\, \alpha^{n+1}_1\nabla p^{n+1} + \gamma_{12}^n(\ve{u}_1^{n+1}-\ve{u}_2^{n+1})
		&=  \rho_1\frac{\alpha_1^{n}}{\tau}\ve{u}_1^{n} +  \rho_1\alpha_1^{n+1}\ve{g}_1^{n+1},\\
	\rho_2\bigg[\frac{\alpha_2^{n+1}+\alpha_2^{n}}{2\tau}\ve{u}_2^{n+1}   + \alpha_2^{n+1}\ve{u}_2^n\cdot\nabla\ve{u}_2^{n+1} + \frac{1}{2}\nabla\cdot(\alpha_2^{n+1}\ve{u}_2^n)\ve{u}_2^{n+1}\bigg]  &-\nabla\cdot(2\mu_2\alpha^{n+1}_2\nabla^{\mathrm{s}}\ve{u}^{n+1}_2)  \\
		+\, \alpha^{n+1}_2\nabla p^{n+1} + \gamma_{12}^n(\ve{u}_2^{n+1}-\ve{u}_1^{n+1})
		&=  \rho_2\frac{\alpha_2^{n}}{\tau}\ve{u}_2^{n} +  \rho_2\alpha_2^{n+1}\ve{g}_2^{n+1},\\	\nabla\cdot(\alpha_1^{n+1}\ve{u}_1^{n+1}+\alpha_2^{n+1}\ve{u}_2^{n+1})&= 0\, ,
	\end{align*}
which is also linearised and unconditionally stable. The comparison in terms of $\alpha_1$ and $p$ is shown in Figure \ref{alpha1}, revealing very good agreement between the fractional-step and monolithic schemes, even though they also differ with respect to the chosen change of variable ($\phi_k$ and $\varphi_k$ for the fractional-step and monolithic cases, respectively).

\begin{figure}[ht!]
		\centering
		\includegraphics[trim = 35 12 95 10,clip, width = 1\textwidth]{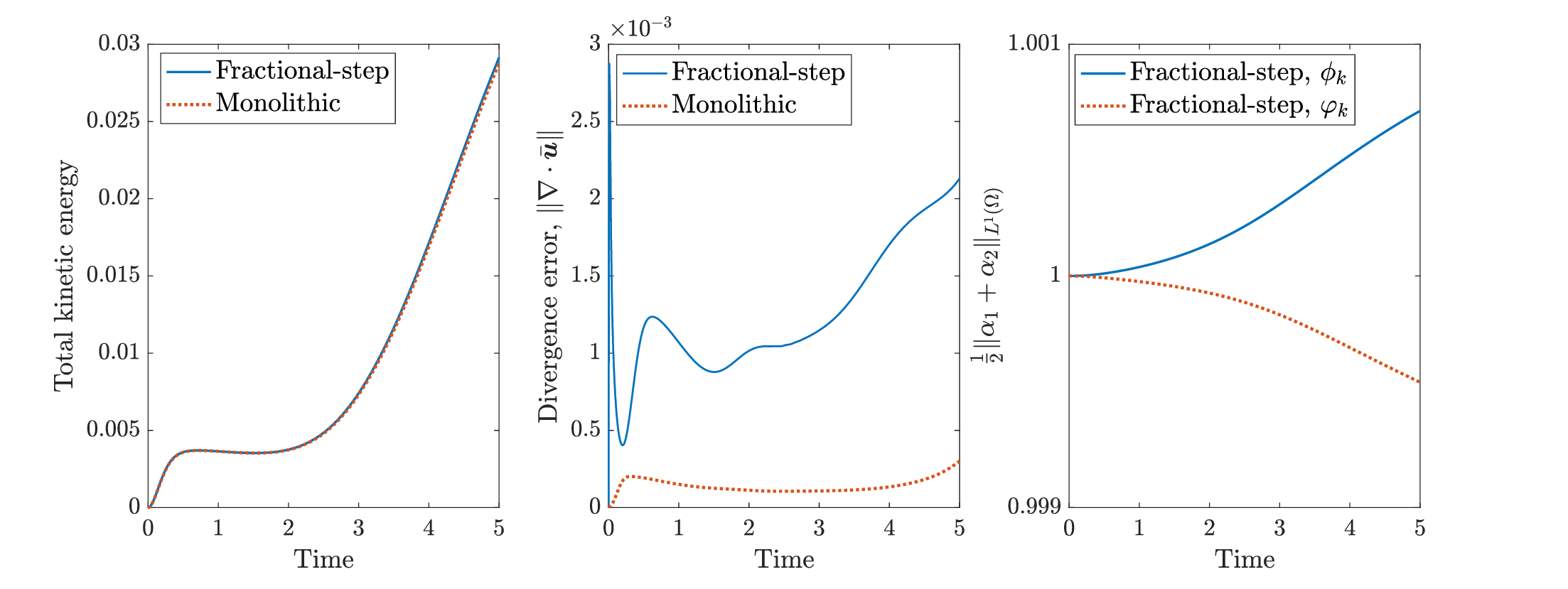}
		\caption{Dispersed Rayleigh--Taylor flow: temporal evolution of three quantities of interest. The results on the right compare two different changes of variable for $\alpha_k$, see Section \ref{sec_convection}.}
		\label{evolution}
	\end{figure}

As a last result, we assess the evolution of three quantities of interest:
\begin{align*}
E_{\mathrm{kinetic}} &= \frac{1}{2}\rho_1\left\|\sqrt{\alpha_1}\,\ve{u}_1\right\|^2 + \frac{1}{2}\rho_2\left\|\sqrt{\alpha_2}\,\ve{u}_2\right\|^2\, ,\\
 \|\nabla\cdot\bar{\ve{u}}\| &= \|\nabla\cdot(\alpha_1\ve{u}_1+\alpha_2\ve{u}_2)\|\, ,\\
 \frac{1}{|\Omega|}\|\alpha_1 + \alpha_2\|_{L^1(\Omega)} &= \frac{\|\varphi_1 \|^2 + \|\varphi_2 \|^2}{|\Omega|}\, .
\end{align*}
While the kinetic energy should change due to gravity and dissipation, the other two quantities should (for the exact solution) remain equal to $0$ and $1$, respectively, throughout time. The results in Figure \ref{evolution} show that the kinetic energy obtained through the fractional-step and monolithic schemes agree well. As expected, the divergence error is smaller in the monolithic case, which enforces $\nabla\cdot\bar{\ve{u}}$ directly (albeit weakly). For the conservation of $\|\alpha_1+\alpha_2\|_{L^1(\Omega)}$, the comparison between the two change-of-variable approaches discussed in Section \ref{sec_convection} confirms the theoretical predictions: both versions are similarly accurate, but the bound-preserving one leads to a (slight) progressive growth over time.

	\section{Concluding remarks}\label{sec_Conclusion}
	This work has presented a simple IMEX splitting scheme for dispersed multiphase flow systems. While each fluid phase obeys its own mass and momentum equations, they are all coupled by the pressure and drag forces. Therefore, decoupling the phases (at each time step) requires treating those forces explicitly. The pressure splitting is a particularly delicate matter, even more so here than in single-phase flows. To achieve that, we have presented a pressure-correction method built upon the incompressibility of the mean (volume-fraction-weighted) velocity field. We have also proposed an explicit treatment of the drag term based on the projected velocities. Moreover, the convective nonlinearities are also eliminated through simple IMEX techniques. Finally, we also showed that one may choose to treat part of the viscous term explicitly to further decouple the spatial components of each phase velocity. The resulting scheme, which we have proved unconditionally stable, is considerably simpler to implement than an implicit one, as it is linearised and fully decoupled---as a matter of fact, only scalar subproblems need to be solved. Although we focus on two-phase applications, our method becomes increasingly more attractive for more fluids because only one pressure problem is solved regardless of the number of phases. Of course, our framework has limitations. For instance, its extension to second order may prove challenging to construct and prove stable. Even first-order convergence may be very difficult to attain in the presence of open boundaries, where standard pressure-correction methods usually fail also for single-phase flows \cite{Guermond2005}. A promising path towards enabling both higher-order convergence and more general boundary conditions is consistent splitting schemes \cite{Pacheco2022CompMech}, which is a topic we are currently working on.

\appendix
\section{Deriving the pressure-correction method}\label{sec_appendix}
Projection schemes can be derived in various ways, and here we follow the pipeline proposed by \citet{BadiaRamon2007} for the homogeneous, single-phase Navier--Stokes equations. We start from the fully implicit Euler discretisation of \eqref{momentum2}--\eqref{divFree}:
\begin{align}
		\frac{\rho_k}{\tau}\left(\alpha_k^{n+1}\ve{U}_k^{n+1}
    -\sqrt{\alpha_k^{n+1}\alpha_k^{n}}\,\ve{U}_k^{n} \right) + \ve{F}^{n+1}_k(\ve{U}_k^{n+1}) + \alpha_k^{n+1}\nabla p^{n+1} &= \rho_k\alpha_k^{n+1}\ve{g}_k^{n+1}\, ,  \label{momentumEuler}\\
		\nabla\cdot\sum_{k=1}^M\alpha_k^{n+1}\ve{U}_k^{n+1} &= 0\, ,\label{divFreeImplicit}
	\end{align}
    where we capitalise $\ve{U}_k^{n+1}$ to distinguish between the implicit and IMEX solutions, and
    \begin{align*}  &\ve{F}^{n+1}_k(\ve{U}_k^{n+1}):=\rho_k\alpha_k^{n+1}\ve{U}_k^{n+1}\cdot\nabla\ve{U}_k^{n+1}  + \frac{\rho_k}{2} \nabla\cdot(\alpha_k^{n+1}\ve{U}_k^{n+1})\ve{U}_k^{n+1}  - \nabla\cdot(2\alpha_k^{n+1}\mu_k\nabla^{\mathrm{s}}\ve{U}_k^{n+1})\\
        &+ \sum_{l=1}^{M}\gamma_{kl}^{n+1}(\ve{U}_k^{n+1}-\ve{U}_l^{n+1}) \, .
    \end{align*}
    Next, we introduce the auxiliary velocities 
    \begin{align}
        \tilde{\ve{U}}_k^{n+1} = \ve{U}_k^{n+1} + \frac{\tau}{\rho_k}\left(\nabla p^{n+1} - \sqrt{\frac{\alpha_k^{n}}{\alpha_k^{n+1}}}\,\nabla p^n\right)
        \label{auxVel}
    \end{align}
and use them in the first term on the left-hand side of \eqref{momentumEuler} to get
    \begin{align}
		\frac{\rho_k}{\tau}\left(\alpha_k^{n+1}\tilde{\ve{U}}_k^{n+1}
    -\sqrt{\alpha_k^{n+1}\alpha_k^{n}}\,\ve{U}_k^{n} \right) + \ve{F}^{n+1}_k(\ve{U}_k^{n+1}) + \sqrt{\alpha_k^{n+1}\alpha_k^{n}}\,\nabla p^{n} &= \rho_k\alpha_k^{n+1}\ve{g}_k^{n+1}\, .  \label{momentumEulerMod}
\end{align}
Now, multiplying \eqref{auxVel} by $\alpha_k^{n+1}$, adding over the phases and taking the divergence yields
\begin{align*}
\nabla\cdot\sum_{k=1}^M\left[\alpha_k^{n+1}\tilde{\ve{U}}_k^{n+1} + \frac{\tau}{\rho_k}\left(\sqrt{\alpha_k^{n+1}\alpha_k^{n}}\,\nabla p^{n} - \alpha_k^{n+1}\nabla p^{n+1}\right)\right]
    &=\nabla\cdot\sum_{k=1}^M\alpha_k^{n+1}\ve{U}_k^{n+1}\\
    &= 0\, ,
\end{align*}
which can be reorganised as a Poisson-like equation for $p^{n+1}$:
\begin{align}
\nabla\cdot\left[\left(\sum_{k=1}^{M}\frac{\alpha_k^{n+1}}{\rho_k}\right)\nabla p^{n+1}\right] 
= \nabla\cdot\left[\left(\sum_{k=1}^{M}\frac{\sqrt{\alpha_k^{n+1}\alpha_k^{n}}}{\rho_k}\right)\nabla p^{n}\right] +\frac{1}{\tau}\nabla\cdot\sum_{k=1}^{M}\alpha_k^{n+1}\tilde{\ve{U}}_k^{n+1}\, .\label{PPEderivation}
\end{align}

Notice that, up to this point, the modification of the implicit system is fully consistent: system \eqref{auxVel}--\eqref{PPEderivation} is equivalent to \eqref{momentumEuler}--\eqref{divFreeImplicit} (all the derivation steps can be reversed). The so-called splitting error comes now, as we introduce the approximation $\ve{F}^{n+1}_k(\ve{U}_k^{n+1}) \approx \ve{F}^{n+1}_k(\tilde{\ve{U}}_k^{n+1})$ and the artificial boundary condition 
\begin{align}
    \left(\sum_{k=1}^{M}\frac{\alpha_k^{n+1}}{\rho_k}\right) \partial_{\vesmall{n}}p^{n+1} = \left(\sum_{k=1}^{M}\frac{\sqrt{\alpha_k^{n+1}\alpha_k^{n}}}{\rho_k}\right)\partial_{\vesmall{n}}p^{n}\ \ \text{on} \ \, \partial\Omega\, ,
\end{align}
which is equivalent to 
\begin{align*}
 \ve{n}\cdot\sum_{k=1}^M(\tilde{\ve{U}}_k^{n+1} -\ve{U}_k^{n+1}) = 0\   \ \text{on} \ \, \partial\Omega\, ,
\end{align*}
due to \eqref{auxVel}. Finally, changing notation from $\ve{U}_k^{n+1}$ to $\hat{\mathbf{u}}_k^{n+1}$ and from $\tilde{\ve{U}}_k^{n+1}$ to $\ve{u}_k^{n+1}$, we have the split system
\begin{align}
    &\frac{\rho_k}{\tau}\left(\alpha_k^{n+1}\ve{u}_k^{n+1}
    -\sqrt{\alpha_k^{n+1}\alpha_k^{n}}\,\hat{\mathbf{u}}_k^{n+1}\right) + \ve{F}^{n+1}_k(\ve{u}_k^{n+1}) + \sqrt{\alpha_k^{n+1}\alpha_k^{n}}\,\nabla p^{n} = \rho_k\alpha_k^{n+1}\ve{g}_k^{n+1}\, ,\\
   & \nabla\cdot\left[\left(\sum_{k=1}^{M}\frac{\alpha_k^{n+1}}{\rho_k}\right)\nabla p^{n+1}\right] 
= \nabla\cdot\left[\left(\sum_{k=1}^{M}\frac{\sqrt{\alpha_k^{n+1}\alpha_k^{n}}}{\rho_k}\right)\nabla p^{n}\right] +\frac{1}{\tau}\nabla\cdot\sum_{k=1}^{M}\alpha_k^{n+1}\ve{u}_k^{n+1}\, ,\\
&\hat{\mathbf{u}}_k^{n+1} = \ve{u}_k^{n+1} + \frac{\tau}{\rho_k}\left(\sqrt{\frac{\alpha_k^n}{\alpha_k^{n+1}}}\,\nabla p^{n}-\nabla p^{n+1}\right),
\end{align}
to be solved sequentially. To obtain the precise version considered in this work, it remains only to introduce the first-order IMEX approximations
\begin{align*}  &\ve{F}^{n+1}_k(\ve{u}_k^{n+1})\approx \rho_k\alpha_k^{n+1}\ve{u}_k^{n}\cdot\nabla\ve{u}_k^{n+1}  + \frac{\rho_k}{2} \nabla\cdot(\alpha_k^{n+1}\ve{u}_k^{n})\ve{u}_k^{n+1}  - \nabla\cdot(2\alpha_k^{n+1}\mu_k\nabla^{\mathrm{s}}\ve{u}_k^{n+1})\\
        &+ \sum_{l=1}^{M}\gamma_{kl}^{n}(\hat{\mathbf{u}}_k^{n}-\hat{\mathbf{u}}_l^{n}) \, ,
    \end{align*}
along with \eqref{Salgado}. This method is formally first-order consistent and has an $\mathcal{O}(\tau^2)$ splitting error, as shown in Section \ref{sec_penalty}

	\section*{Acknowledgments}
	\noindent DRQP acknowledges funding by the Federal Ministry of Education and Research (BMBF) and the Ministry of Culture and Science of the German State of North Rhine-Westphalia (MKW) under the Excellence Strategy of the Federal Government and the Länder.
	
	\bibliographystyle{unsrtnat} 
	\bibliography{references}%
	
\end{document}